\magnification=\magstep1    
\input amssym.def
\input amssym.tex
\hsize=5.9truein
\overfullrule=0pt
\def\Q{{\Bbb Q}}
\def\Z{{\Bbb Z}}
\def\Br{\mathop{\rm Br}\nolimits}
\def\Gal{\mathop{\rm Gal}\nolimits}
\def\Hom{\mathop{\rm Hom}\nolimits}
\def\Pic{\mathop{\rm Pic}\nolimits}
\def\ram{\mathop{\rm ram}\nolimits}
\def\Spec{\mathop{\rm Spec}\nolimits}
\def\diagram{\def\normalbaselines{\baselineskip21pt\lineskip1pt
	\lineskiplimit0pt}}
\def\proof{\noindent{\it Proof.}\quad}
\def\blackbox{\hbox{\vrule width6pt height7pt depth1pt}} 
\def\qed{~\hfill\blackbox\medskip}
\outer\def\Demo #1. #2\par{\medbreak\noindent {\it#1.\enspace}
	{\rm#2}\par\ifdim\lastskip<\medskipamount\removelastskip
	\penalty55\medskip\fi}
\def\hangbox to #1 #2{\vskip1pt\hangindent #1\noindent \hbox to #1{#2}$\!\!$}

\pageno=0 
\footline{\ifnum\pageno=0\hfill\else\hss\tenrm\folio\hss\fi}
\topinsert\vskip1.8truecm\endinsert
\centerline{\bf Cyclic algebras over $p$-adic curves}
\vskip6pt
$${\vbox{\halign{\hfil\hbox{#}\hfil\qquad&\hfil\hbox{#}\hfil\cr
$$David J. Saltman$^*$\cr
Department of Mathematics\cr
The University of Texas\cr
Austin, Texas 78712\cr}}}$$
\vskip16pt
\vskip3pt
{\narrower\smallskip\noindent
{\bf Abstract} 
In this paper we study division algebras over the function fields 
of curves over $\Q_p$. The first and main tool is to view 
these fields as function fields over nonsingular $S$ 
which are projective of relative dimension 1 over 
the $p$ adic ring $\Z_p$. A previous paper showed such division algebras 
had index bounded by $n^2$ assuming the exponent was $n$ and $n$ was prime 
to $p$. In this paper we consider algebras of degree 
(and hence exponent) $q \not= p$ and show these 
algebras are cyclic. We also find a geometric criterion 
for a Brauer class to have index $q$.  

\bigskip

\medskip

\smallskip}

\footnote{}{*The author is grateful for support under NSF grants DMS-9970213 
and DMS-0401468.}
\vfill\eject

\leftline{Introduction}
\medskip 
In [S], this author studied division algebras over the following 
fields. Let $K$ be a field finite over $\Q_p(t)$, for the $p$ adic field 
$\Q_p$. That is, suppose there is a $p$ adic field $K'$ and a curve 
$C$ defined over $K'$ such that $K = K'(C)$. Let $n$ be prime to $p$. 
In [S] we studied division algebras $D/K$ 
(meaning $K$ is the center of $D$) and showed 
that if their order in the Brauer group was $n$, then their degree 
was no more than $n^2$.  

This paper is motivated by the idea that there are 
further interesting things to say about division 
algebras over these fields $K$. For example, 
suppose $D/K$ has degree $q^2$, for a prime $q \not= p$, 
and order $q$ in the Brauer group. The techniques of [S] 
show that, assuming $K$ has a primitive $q$ root of one, 
then $D/K$ is an abelian crossed product (e.g. [LN] p. 37). 
For this and other more obvious reasons, it is of 
interest to study $D/K$ of prime degree $q$. 
The important question is whether these $D$ are cyclic 
algebras, and the answer we provide here is that  
such $D$ are cyclic, whether or not there are $q$ roots of one (5.1).

The first important step here, as in [S], is to observe 
that $K$ is the function field of a regular surface $S$ projective 
over $\Spec(\Z_p)$, where $\Z_p$ is the ring of $p$-adic integers. 
Thus much of this paper will have a geometric character, 
as the geometry imposed on $S$ by $D$ needs to be explicated and 
understood. Let me apologize in advance to geometers 
for the proofs I may provide for well known facts. 
Part of the intended audience of this paper consists of people 
primarily concerned with division algebras. I have chosen, 
therefore, to provide proofs of any facts that cannot be found 
in standard texts like Hartshorne and EGA. 

Let me review briefly the structure of the paper. Our approach will 
be to prove as much as we can about Brauer classes over surfaces, 
and only use the strong condition on $S$ above when needed at the end. 
In more detail, in this introductory section we review some material about Brauer 
groups, ramification, cyclic extensions, etc. 
Section one is a general geometry section. In the first half 
we review facts about surfaces $S$ projective over $\Z_p$, 
and in the second half we consider a cohomology group 
$H^1(X,{\cal O}^*_P)$ over a much more general scheme $X$. 
The point is to do ``divisor theory'' while controlling 
behavior at finitely many points. In section two we assume 
the ground field has a primitive $q$ root of one, and study 
the geometry of the ramification of a Brauer class of order $q$. 
In section three we remove the assumption on roots of unity. 
In section four we consider the behavior of ``residual'' 
classes, and in section five we prove the main results.

Let $q$ be a fixed prime unequal to $p$ throughout this paper. 
Let $\mu_q$ be the group of $q$ roots of one over any field. 
We denote by $G_F$ the absolute Galois group of a field $F$. 
That is, $G_F$ is the Galois group of $F$ in its separable closure.
If $\mu_q \subset F^*$, there is a pairing $G_F \times F^* \to \mu_q$ 
defined by sending $(\sigma, u) \to \sigma(u^{1/q})/u^{1/q}$. 
If $F$ is a finite field containing $\mu_q$, 
then the Frobenius defines a canonical 
generator of $G_F$ and so the Frobenius 
defines a homomorphism $Fr: F^* \to \mu_q$. 

Recall that if $K$ is a field, 
the Brauer group $\Br(K)$ consists of equivalence classes $[A]$
of central simple algebras $A/K$, and each such class contains a unique 
division algebra. If $\alpha \in \Br(K)$, then the order 
of $\alpha$ is its order in the Brauer group, and the 
index of $\alpha$ is the degree (i.e. square root of the 
dimension) of the associated division algebra over $K$. 
A cyclic algebra is a central simple algebra $A/K$ of degree 
$n$ containing $L$ where $L/K$ is cyclic Galois of degree 
$n$ ($L$ need not be a field). All cyclic algebras have the form 
$A = \Delta(L/K,\sigma,a)$ where $L/K$ is cyclic Galois, 
$\sigma \in \Gal(L/K)$ is a generator, and $a \in K^*$ (e.g. [LN] p. 49).  
Note that $\Delta(L/K,\sigma,a) \cong \Delta(L/K,\sigma^s,a^s)$ 
where $s$ is prime to the degree of $L/K$.  
If $K' \supset K$ is a field 
extension, recall that $\alpha \in \Br(K)$ is split by 
$K'$ if it is in the kernel of the natural map 
$\Br(K) \to \Br(K')$ given by $[A] \to [A \otimes_K K']$. 
Perhaps the most important fact about cyclic algebras 
we need is the well known theorem of Albert: 

\proclaim Proposition 0.1. Suppose $A/K$ is a central simple algebra 
of prime degree $q$. Then $A$ is a cyclic algebra if and only if 
there is a $\pi \in K^*$ such that $K' = K(\pi^{1/q})$ 
splits $[A]$. 

\proof 
The description of cyclic algebras above shows that they contain 
such Kummer maximal subfields, and such a subfield necessarily splits 
$A$. 
Thus the ``only if'' part is done. If such a $K'$ splits $A$, 
by [LN] p. 25 it is isomorphic to a subfield of $A$. This result 
now follows from [A] p. 77.~\qed  

When $F$ contains $\rho$, a generator of $\mu_q$, 
all cyclic algebras over $F$ have the following 
form. If $a,b \in F^*$ then one can define the symbol algebra 
$(a,b)_{q,F,\rho}$ as the central simple $F$ algebra generated by 
$x$, $y$ satisfying the relations $x^q = a$, $y^q = b$ and 
$yx = \rho{xy}$. Just as with general cyclic algebras, 
we have that $(a,b)_{q,F,\rho} \cong (a,b^s)_{q,F,\rho^s}$ 
where $s$ is prime to $q$. We will often drop all or a subset 
of the $q,F,\rho$ subscript 
because $q$ is fixed throughout the paper, $F$ is usually clear, 
and $\rho$ is often fixed in advance. 
We will also write $(a,b) \in \Br(F)$ for the Brauer group 
element represented by the algebra $(a,b)$ (and called a  {\bf symbol class}). 

If $R$ 
is a discrete valuation domain with field of fractions 
$q(S) = K$ and residue field $F$ of characteristic $p$, 
then there is the well known ramification map (e.g. [Se] p. 186) 
$$\ram: \Br(K)' \to \Hom(G_F,\Q/\Z)'$$ 
where for a torsion abelian group $A$, $A'$ refers to the 
prime $p$ part of $A$. 
Note that any $q$ order element $\phi \in \Hom(G_F,\Q/\Z)$ 
can be represented by a pair $L/F,\sigma$ where the kernel of 
$\phi$ has fixed field $L$ and $\sigma$ is the generator of 
$C_q = \Gal(L/F)$ which maps to $1/q + \Z$. In this paper 
ramification will be frequently written this way. 

This ramification 
map is almost completely determined by the following 
two observations. First, let $\tilde K$ be the completion 
of $K$ with respect to $R$. Then the ramification map factors 
as $\Br(K)' \to \Br(\tilde K)' \to \Hom(G_F,\Q/\Z)'$ 
where the first map is the usual restriction on Brauer groups 
and the second map is the ramification associated to 
the valuation on $\tilde K$. Second, assume $K = \tilde K$ 
is complete. Suppose $L/K$ is cyclic unramified 
of degree prime to $p$, with generator $\sigma \in \Gal(L/K)$. 
Let $\bar L/F,\bar \sigma$ be the residue extension and corresponding 
generator. Then the ramification of the cyclic algebra 
$\Delta(L/K,\sigma,\pi)$ is $\bar L/F,\bar \sigma$ when $\pi$ 
is any prime element of $K$.   

In particular, 
assume $F$ contains $\mu_q$ and we fix a generator $\rho \in \mu_q$.  
Then the $q$ torsion part of $\Hom(G_F,\Q/\Z)$ can be identified with 
$F^*/(F^*)^q$. In detail, the pair $L/F, \sigma$ is identified 
with $a(F^*)^q$ where $L = F(a^{1/q})$ and $\sigma(a^{1/q})/a^{1/q} = 
\rho$.  
Thus a $q$ torsion element of $\Hom(G_F,\Q/\Z)$ 
will sometimes represented by an element of $F^*/(F^*)^q$ 
or by $a^{1/q}$ for some $a \in F^*$. 
With all of this, there is an easy way to write the 
ramification of a symbol class $(a,b)$. The following 
result is well known and computable from the above or, 
for example, [LN] p. 68. 

\proclaim Lemma 0.2. Suppose $R \subset K$, $F$ 
are as above, $\rho \in \mu_q \subset K$ is fixed, and $(a,b) \in \Br(K)$  
is a symbol class. Let $d: K^* \to \Z$ be the valuation associated to $R$. 
Then $\ram((a,b)) = (\bar u)^{1/q}$ where 
$u = (-1)^{d(a)d(b)}a^{d(b)}/b^{d(a)}$ 
and where $\bar u$ refers to the image of $u$ in $F^*$.  

Suppose we have a field $K$ which is the function 
field of a normal integral scheme $X$ of finite type over a Noetherian 
ring. 
Let $\alpha \in \Br(K)$. For each irreducible divisor 
$D \subset X$ let $R_D$ be the stalk of the structure sheaf 
of $X$, which is a discrete 
valuation domain. There are only finitely many 
$D_i$ where $\alpha$ has nontrivial ramification $L_i/F(D_i),\sigma_i$. 
The set of $D_i$ where $\alpha$ is ramified 
is called the {\bf ramification locus} of $\alpha$. 
The set of $D_i$ paired with the ramification 
$L_i/F(D_i),\sigma_i$ of $\alpha$ 
at each $D_i$ is called the {\bf ramification data} of $\alpha$. 

Much of this paper is about splitting ramification so it is important 
we describe how this is done. Let $R \subset K$ be a discrete 
valuation domain of $K$ (meaning $K$ is the field of fractions 
of $R$) and let $F$ be the residue field of $R$. 
Let $L/K$ be a finite separable extension 
field and let $\{S_i\}$ be the 
(necessarily finite) set of discrete valuation domains 
of $L$ which extend $R$. Let $F_i$  be the residue field 
of $S_i$ and $e_i = e(S_i/R)$ the ramification index. 
Let $\ram_i: \Br(L)' \to \Hom(G_{F_i},\Q/\Z)'$ and 
$\ram: \Br(K)' \to \Hom(G_F,\Q/\Z)'$ be the respective 
ramification maps. 

\proclaim Lemma 0.3. The following diagram commutes:
$$\diagram 
\matrix{
\Br(L)'&\buildrel{\sum \ram_i}\over\longrightarrow&\oplus \Hom(G_{F_i},\Q/\Z)'\cr
\llap{$\iota$}\uparrow&&\uparrow\rlap{$\sum e_i$}\cr
\Br(K)&\buildrel{\ram}\over\longrightarrow&\Hom(G_F,\Q/\Z)'}$$
where $\iota$ is the restriction and $e_i: \Hom(G_F,\Q/\Z)' 
\to \Hom(G_{F_i},\Q/\Z)'$ is the natural map multiplied by 
the integer $e_i$. 

If $\iota(\alpha)$ is unramified at all $S_i$ we say {\bf L/K 
splits all the ramification of $\alpha$ at R}. 
We are particularly interested in the case $L/K$ above is a 
field extension of prime degree $q$ unequal to the residue 
characteristic. 

\proclaim Corollary 0.4. Let $L/K$ in 0.3 be  
of prime degree $q$ unequal to the residue characteristic. 
Assume $\alpha \in \Br(K)$ has ramification $F'/F,\sigma$ 
of order $q$. Then $\iota(\alpha)$ is unramified at all the 
$S_i$ if and only if there is a unique extension, $S$, 
of $R$ to $L$ and one of the following two exclusive 
conditions hold:
\smallskip
i)  $L/K$ is totally and tamely ramified. 
\smallskip
ii) $L/K$ is unramified 
at $R$ and the residue field of $S$ is $F'$. 

\proof 
If $S$ is not unique, all ramification degrees and all 
residue extension degrees are prime to $q$ 
and $L$ cannot split the ramification at any extension. 
On the other hand, suppose $R$ extends to a unique $S$. 
If $L/K$ is ramified, $q = e(S/R)$, and $\iota(\alpha)$ 
has zero ramification at $S$. If $S/R$ is unramified, let 
$F''$ be the residue field of $S$, so $F''/F$ is of degree 
$q$. Then $F'' \supset F'$ if and only if $F'' = F'$ if and only 
if $\iota(\alpha)$ has zero ramification at $S$.~\qed

If $L/K$ as in 0.4 satisfies i) we say it splits $\alpha$ 
{\bf by ramification} 
and if $L/K$ satisfies ii) we say it splits 
$\alpha$ {\bf by residues}. 

Let us make one more definition. 
Suppose $\alpha \in \Br(K)$, $K$ has a discrete valuation 
$R$, and $L/K$ splits the ramification of $\alpha$ 
at $R$ and is totally ramified, which includes that $R$ extends uniquely. 
If $S$ is that unique extension, and $\alpha_L = \alpha \otimes_K L$ 
is the image of $\alpha$ in $\Br(L)$, then 0.3 shows that 
$\alpha_L \in \Br(S)$. If $F$ is the residue field of 
$S$ and hence of $R$, then $\alpha_L$ has an image 
$\beta_R \in \Br(F)$ we call the {\bf residual Brauer class} 
of $\alpha$ at $R$ with respect to $L$. 

We can make the following observation about $\beta_R$. 

\proclaim Proposition 0.5. Suppose $\alpha$, $R$, $K$ and $L$ 
are as above, and let $F'/F,\sigma$ be the nonzero ramification 
of $\alpha$ at $R$. Assume $L/K$ has degree $q$. Suppose $\alpha$ has index $q$, 
meaning it is represented by a division algebra of degree $q$. 
Then the residual Brauer class $\beta_R$, with respect to any $L$, 
is split by $F'$.  

\proof At the completion $\alpha$ must still have index $q$, 
so it suffices to prove this under the assumption that 
$K$ is complete with respect to $R$. In addition, 
it suffices to show this after we adjoin a $q$ root of one. 
Thus we may assume $K$ contains a primitive $q$ root 
of one. Since $L/K$ is totally and tamely ramified, 
it is cyclic of degree $q$. But then $\alpha = \alpha' + (K'/K,\pi)_{q,K}$ 
where $K'/K$ is the unramified extension with residue extension 
$F'/F$, $L = K(\pi^{1/q})$ and $\alpha' \in \Br(R)$ with image $\beta_R$. 
By e.g. [JW] p. 161, if $F'$ does not split $\beta_R$ then $\alpha$ 
has index bigger than $q$.~\qed 

In the rest of this paper, the  $R$ of 0.5 will sometimes be defined 
by a curve $C$ on a surface $S$, and and in that case we will write the 
residual Brauer class as $\beta_C$. 

It will later be important to determine how this residual class 
$\beta_R$ depends on the choice of $L$. To this end, 
let $R \subset K$ be a discrete valuation domain with field 
of fractions $K$ as above. 

\proclaim Proposition 0.6. 
Suppose $\alpha \in \Br(K)$ of order $q$ has ramification $F'/F,\sigma$ 
at $R$. Let $L = K(\pi^{1/q})$ 
for $\pi$ a prime of $R$. Also set $L' = K((u\pi)^{1/q})$ 
where $u$ is a unit of $R$. Let $\beta_R$, $\beta_R'$ 
be the respective residual classes of $\alpha$ defined by
$L$ and $L'$. Then $\beta_R' = \beta_R + \Delta(F'/F,\sigma,\bar u^{-1})$ 
where $\bar u$ is the image of $u$ in $F^*$.

\proof 
Just as above, 
to prove this we can assume $K$ is complete with respect to $R$. 
Let $K'/K$ be unramified 
with residue extension $F'/F$. Then 
$\alpha = \alpha' + \Delta(K'/K,\sigma,\pi) = 
\alpha' + \Delta(K'/K,\sigma,u^{-1}) + \Delta(K'/K,\sigma,u\pi)$. 
Since $\beta_R$ is the image of $\alpha'$ and $\beta_R'$ 
is the image of $\alpha' + \Delta(K'/K,\sigma,u^{-1})$, 
we are done.~\qed 

Sometimes it will not be convenient to have $L$ written as $K(\pi^{1/q})$ 
with $\pi$ a prime but only with $\pi$ having prime to $q$ valuation. 
The following is obvious. 

\proclaim Corollary 0.7. Let $\alpha \in \Br(K)$ and $R$, $F$, 
$F'/F,\sigma$  
be as above. Suppose $v$ is the valuation of $R$ and $\pi \in K$ 
satisfies $v(\pi) = s$ which is prime to $q$. Set 
$L = K(\pi^{1/q})$ and let $\beta_R$ be the corresponding 
residual Brauer class. Suppose $u \in R^*$ has image and $L' = 
K((u\pi)^{1/q})$. If $\beta_R'$ is the residual class 
with respect to $L'$, then 
$\beta_R' = \beta_C + \Delta(F'/F,\sigma,\bar u^{-t})$ 
where $st - 1$ is divisible by $q$. 

\Demo Remark. Suppose $\alpha \in \Br(K)$, $R$, and $F'/F$ are as in 
0.7. If $\alpha$ has index $q$, we know by 0.5 
that $F'$ splits $\beta_{R}$. The converse is false but 
0.7 makes the following clear. If for one choice of $L$, 
$\beta_R$ is split by $F'$, then this is true for all 
choices of $L$. When this happens, we say 
the {\bf residual classes of $\alpha$ are split by the ramification}. 
  
Suppose next that $K = k(C)$ is the function field of a curve over 
a finite field $k$. Then the set of discrete valuations on $K$ is 
exactly the set of points on $C$. If $R$ is any such discrete 
valuation, then the residue field $R/P$ is a finite field 
and hence $\Hom(G_{R/P},\Q/\Z)$ can be identified with 
$\Q/\Z$ using evaluation on the Frobenius. Thus there is a map 
$\Br(K) \to \oplus_{P \in C} \Q/\Z$. Note that since all finite 
fields are perfect, we can define this map even on the $p$ 
primary part of $\Br(K)$ (e.g. [Se] p. 186). 
From class field theory we know (e.g. [R] p. 277): 

\proclaim Theorem 0.8. There is an exact sequence 
$$0 \to \Br(K) \to \oplus_{P \in C} \Q/\Z \to \Q/\Z \to 0$$
where $\oplus_{P \in C} \Q/\Z \to \Q/\Z$ is the 
summation map. 

If $\alpha \in \Br(K)$ then the image of $\alpha$ in the copy 
of $\Q/\Z$ corresponding to $P \in C$ we call the {\bf residue} 
of $\alpha$ at $P$. 

One consequence of 0.8 is that, over $K = k(C)$ as in 0.8, 
splitting all ramification 
is equivalent to splitting. This is false in general, 
though of course splitting implies splitting all ramification. 
The most important fact about the fields $K$ that are the focus 
of this paper is that for them also, 
splitting all ramification implies splitting. 

\proclaim Theorem 0.9. Suppose $S$ is a surface projective and 
regular over $\Spec(\Z_p)$. Let $K$ be the function field of 
$S$. If $\alpha \in \Br(K)$ has trivial ramification 
at all discrete valuations lying over $\Z_p$, then 
$\alpha = 0$. 

\proof It suffices to show $\Br(S) = 0$. By [G] p. 98,  
$\Br(S) \cong \Br(\bar S)$ and so it suffices to show 
$\Br(\bar S) = 0$. By, for example, the argument of 
[S] p. 40 it suffices to show $\Br(C) = 0$ for $C$ 
any complete nonsingular curve over a finite field, 
and this is a part of 0.8.~\qed 

Having discussed ramification of algebras, let us consider 
that of cyclic extensions. Let $R$ be a discrete valuation 
domain with residue field $F = R/{\cal M}$ and field of fractions 
$K = q(R)$. Suppose $L/K$ is a cyclic Galois extension 
of prime degree $q$ with generator $\sigma$ of its Galois 
group. We assume $q$ is not the characteristic of $k$ 
and $\mu_q$ is the group of $q$ roots of one over $K$. 
We need to define the ramification $\rho \in \mu_q$
of $L/K,\sigma$ at $R$. 
If $L/K$ is unramified, of course $\rho = 1$. If $L/K$ 
is ramified, let $\tilde K$ be the completion of $K$ 
with respect to $R$ and $\tilde L = L \otimes_K \tilde K$. 
Then $\tilde L$ is a field. Since $\tilde L/\tilde K$ 
is a totally and tamely ramified extension, it follows that 
$\mu_q \subset \tilde K$ and hence $\mu_q \subset k$. Furthermore, 
$\tilde L/\tilde K, \sigma$ has the form $\tilde K((\pi)^{1/q})$ 
for some prime element $\pi$. Note that $\pi^{1/q} \in \tilde L$ 
is a prime element of $\tilde L$. We set the ramification 
$\rho = \sigma(\pi^{1/q})/\pi^{1/q}$ viewed as a root 
of unity over $k$. It is useful to note that this $\rho$ can be defined 
using any prime element of $\tilde L$ and hence of $L$. In fact, 
suppose $\delta$ is a prime element of $\tilde L$. Then 
$\delta = u\pi^{1/q}$ for a unit $u$ of $\tilde L$. 
Since $\sigma$ acts trivially on the residue field 
of $\tilde L$, it follows that $\rho$ is the image of 
$\sigma(\delta)/\delta$ in the residue field of $L$. 

The ramification of a cyclic extension can be used to express 
the ramification of a cyclic algebra as follows.
Suppose $K$ is a field with a discrete valuation domain 
$R$ and $\alpha = \Delta(L/K,\sigma,u)$ is of degree $q$ where $u \in R^*$ 
and $u$ has image $\bar u$ in the residue field $F$ of $R$. 
If $L/K$ is unramified then $\alpha$ has 0 ramification. 
If not, $F$ contains a 
primitive $q$ root of one. Let $\rho$ be the ramification 
of $L/K,\sigma$ at $R$. The following is easy. 

\proclaim Lemma 0.10.  The ramification of $\alpha$ 
is described by $F(\bar u^{1/q}),\sigma'$ where \break
$\sigma'(\bar u^{-1/q})/\bar u^{-1/q}) = \rho$, and 
$\rho$ is the ramification of $L/K,\sigma$ at $R$.

In a couple of places in this paper we will need to know 
certain discrete valuations exist, beyond those that 
arise from blowing up points. To this end, 
let $R$ be a local domain with field of fractions $K = q(R)$. 
A discrete valuation $d: K^* \to \Z$ of $K$ is said to {\bf lie over} 
$R$ if $d(R) \geq 0$ and $d(R) \not= \{0\}$. 
If $P = \{r \in R | d(r) > 0\}$ then $P$ 
is a nonzero prime and we say $d$ lies over $P$. 
If $R$ is a domain and $L/K$ splits all the ramification 
at any discrete valuation lying over $R$ we say $L/K$ 
splits all the ramification of $\alpha$ at $R$. 

\proclaim Lemma 0.12. Suppose $R$ is a two dimensional 
local regular domain with parameters 
$\pi,\delta$, residue field $k = R/M$,  
and field of fractions $K$. Let $T$ be transcendental 
over $K$. Suppose $a,b \in \Z$ are positive 
integers. Then there is a valuation $d: K(T)^* \to \Z$ on $K(T)$ 
with the following properties. 
First of all, $d(T) = 1$, $d(\pi) = a$ and $d(\delta) = b$. 
Secondly, the residue field of $d$ is $k(\pi',\delta')$ 
where $\pi'$ is the image of $\pi/T^a$, $\delta'$ is the image of 
$\delta/T^b$, and $\pi'$, $\delta'$ are transcendental over $k$. 

\proof Form the polynomial ring $R[T,\pi'',\delta'']$. 
Let $R'$ be the localization of this polynomial 
ring at the maximal ideal generated by $\pi,\delta,T,\pi'',\delta''$, 
so $R'$ is also a regular local domain. 
Then $T^a\pi'' - \pi, T^b\delta''-\delta, T, \pi'', \delta''$ 
clearly generate the maximal ideal of $R'$ and hence form an $R$ 
sequence. Let $R_1 = R'/(T^a\pi'' - \pi, T^b\delta'' - \delta)$ 
which is a regular local ring with parameters we can identify with 
$T$, $\pi''$, $\delta''$. 
Then $R \subset R_1$ and $R_1$ has field of fractions 
$K(T)$. Let $S$ be the discrete valuation ring formed 
by localizing $R_1$ at its prime $T$ and let $d$ be the associated 
valuation. Clearly $d(\pi) = a$, $d(\delta) = b$, and the residue 
field of $S$ is $k(\pi', \delta')$ where 
$\pi'$, $\delta'$ are the images of $\pi''$, $\delta''$ 
and are transcendental over $k$.~\qed 

We will make frequent use of the well known fact (e.g. [E] p. 487) 
that a regular local ring is a UFD. In fact, we will 
need a very slight generalization: 

\proclaim Lemma 0.12. Suppose $R$ is a regular semilocal 
ring. Then $R$ is a UFD. 

\proof It suffices to show that every height one prime $P$ 
is principal. But if $M \subset R$ is a maximal ideal, 
$PR_M$ is a height one prime and hence principal. That is, 
$P$ is locally free, therefore projective, and therefore 
free of rank one since $R$ is semilocal.~\qed 

\leftline{Section One: The surface} 
\medskip
Let $S \to \Spec(\Z_p)$ be projective, regular, excellent, 
flat of relative dimension one. Let $\bar S$ 
be the set theoretic inverse image of the closed point of $\Spec(\Z_p)$ 
with the reduced induced structure. We also assume $\bar S$ has nonsingular 
components and only normal crossings. In this section we review some 
general facts about this situation, which we will apply to the 
Brauer group in subsequent sections. 

First of all let us consider closed points on $S$, by which we mean 
codimension 2 closed points. It is easy to see that 
all such points lie on $\bar S$. Next, we consider codimension 
1 points which we call curves. $\bar S$ is the finite union 
of curves. If $E \subset S$ is any other curve, it lies 
over the generic point of $\Z_p$ and thus defines a  
point of the $\Q_p$ curve $S \times_{\Z_p} \Q_p$. 
The restriction $E \to \Spec{\Z_p}$ is surjective, projective, 
of relative dimension 0 and so must be finite. 
Thus ([H] p. 280) $E$ is affine with affine ring, $R$, a domain 
finite over $\Z_p$. The Henselian property of $\Z_p$ 
shows that $R$ has 0 and one other prime ideal which lies 
over $p\Z_p$. That is, $E$ has a generic 
point and exactly one closed point. 
We call such $E$ geometric curves of $S$.

We observe and recall the well known fact that points 
of $\bar S$ lift nicely to $S$. 

\proclaim Lemma 1.1. a) Let $P \in \bar S$ be a 
(nonsingular) point on a single 
component. There is a nonsingular geometric 
curve $E \subset S$ such that $P$ is the multiplicity one 
intersection of $E$ and $\bar S$. 
\smallskip
b) If $P \in \bar S$ 
is a point on two components, there is a nonsingular geometric 
$E$ which meets each component with multiplicity one 
at $P$. 

\proof 
Let $R = {\cal O}_{S,P}$ be the stalk at $P$ and $M_P$ the maximal ideal. 
In a), 
$p = \delta^ru$ where $u \in R^*$ and $\delta$ is a prime 
of $R$. There is an $x \in R$  such that $(\delta,x) = M_P$. 
Then $R/(x)$ is a DVR, contains 
$\Z_p$, and so must be the integral closure of $\Z_p$ in the field 
of fractions of $R/(x)$. In particular, $R/(x)$ is finite over $\Z_p$. 
If $\Spec(R') \subset S$ is an affine open containing $P$, 
then $R' \subset R$ and $R$ is a localization of $R'$. 
$R'/((x) \cap R') \subset R/(x)$ is also finite over $\Z_p$ 
and so has a unique maximal ideal. The extension 
$R'/((x) \cap R') \subset R/(x)$ is localization at that maximal ideal 
and so $R'/((x) \cap R') = R/(x)$ and this ring represents a nonsingular curve 
geometric curve $E$ 
in $S$ with multiplicity one intersection with $\bar S$ at $P$. 
Since $E$ has a single closed point, this is the only place it intersects 
$\bar S$. 

In b), $p = \delta^r\delta'^su$ where $u \in R^*$ 
and $(\delta,\delta')$ is the maximal ideal of $R$. 
We can now choose $x = \delta + \delta'$ and proceed as above.~\qed

The next issue to concern us is the relation of $\Pic(S)$ 
and $\Pic(\bar S)$. There is a natural map $\Pic(S) \to 
\Pic(\bar S)$ which cannot be an isomorphism but is close enough 
for our needs. 

\proclaim Theorem 1.2. Suppose $\pi: S \to \Spec(\Z_p)$ 
and $\iota: \bar S \to S$ are as above.  
Let $m$ be an integer prime to $p$.  
Then the induced map $\Pic(S) \to \Pic(\bar S)$ 
is a surjection and induces an isomorphism $\Pic(S)/m\Pic(S) 
\cong \Pic(\bar S)/m\Pic(\bar S)$. 

We begin the proof with a proposition. 

\proclaim Proposition 1.3. 
\smallskip 
i) Let $X$ be a scheme and ${\cal J} \subset {\cal O}_X$ 
an ideal sheaf. Let $Y \to X$ be the closed subscheme 
defined by ${\cal J}$. Suppose ${\cal F}$ is a coherent sheaf on 
$X$ with ${\cal J}{\cal F} = 0$. Then $H^i(X,{\cal F}) = H^i(Y,{\cal F})$. 
It is also true that $H^i(X,({\cal O}_X/{\cal J})^*) = H^i(Y,{\cal O}_Y^*)$.  
\smallskip 
ii) Let $X$ be a scheme and ${\cal J} \subset {\cal O}_X$ 
a nilpotent ideal sheaf. Let $Y \to X$ be the closed subscheme 
defined by ${\cal J}$. Assume $Y$ has dimension one, and that the 
integer $m$ is invertible in ${\cal O}_X$. 
Then $\Pic(X) \to \Pic(Y)$ is surjective and induces 
an isomorphism $\Pic(X)/m\Pic(X) \cong \Pic(Y)/m\Pic(Y)$. 

\proof To prove i), note that 
$f_*({\cal F}) = {\cal F}$ and $Y \to X$ is affine 
so exer. 8.2 p. 252 of [H] shows this. The last sentence 
of i) follows similarly.

Turning to ii), by induction we may assume ${\cal J}^2 = 0$. 
There is an exact sequence of abelian group sheaves on $X$: 

$$1 \to {\cal J} \to {\cal O}_X^* \to ({\cal O}_X/J)^* \to 1.$$ 

By i) and [H] p. 208, $H^2(Y, {\cal J}) = 0$ and 
$\Pic(Y) = H^1(X,({\cal O}_X/{\cal J})^*)$. 
It follows from the long exact sequence and i) that 
$\Pic(X) \to \Pic(Y)$ is surjective. 
Also $H^1(X,{\cal J})$ is a module over the ring of global sections 
of $X$, implying that multiplication by $m$ is an isomorphism. 
If $\alpha \in \Pic(X)$ maps to $m\Pic(Y)$, then by the 
surjectivity, there is a $\alpha'$ such that $\alpha - m\alpha'$ 
is the image of $\beta \in H^1(X,{\cal J})$. Since $\beta = m\beta'$ 
for a unique $\beta'$, we have $\alpha = m\alpha' + m\beta''$ 
where $\beta''$ is the image of $\beta'$.~\qed 

We now turn to the proof of 1.2. 

\proof By 1.3, we can replace $\bar S$ with $S_1 \subset S$, 
the subscheme defined by $p{\cal O}_S$. Let $S_n$ be the subscheme 
defined by $p^n{\cal O}_S$. Let ${\cal I}_n \in \Pic(S_n)$ be a previously 
defined line bundle. By 1.3, there is a line bundle 
${\cal I}_{n+1}$ on $S_{n+1}$ such that 
${\cal I}_{n+1}/p^n{\cal I}_{n+1} = {\cal I}_n$. 
By the Grothendieck existence theorem ([EGA] III 5.1.6) 
there is a line bundle ${\cal J}$ on $S$ with 
${\cal J}/p^n{\cal J} = {\cal I}_n$. 

There is another way to view this surjectivity result. 
Since $\bar S$ is a union of smooth curves with 
normal crossings, an element of $\Pic(\bar S)$ can be represented 
as a Cartier Divisor and hence as a sum of points on these curves 
that avoid the intersection points (use 1.5 without circularity). 
Let $P$ be one of these points. Choose $E$ as in 
1.1. Then $E$ defines a divisor, and hence an element 
of $\Pic(S)$ which is the preimage of the element of 
$\Pic(\bar S)$ corresponding to $P$. 

Next, we turn to the injectivity modulo $m$ powers. 
Suppose ${\cal J} \in \Pic(S)$ maps to ${\cal I}^m \in \Pic(\bar S)$. 
Then by lifting ${\cal I}$ we may assume ${\cal J}$ maps to the identity in 
$\Pic(\bar S)$. That is, it suffices to show that the 
kernel of $\Pic(S) \to \Pic(\bar S)$ is $m$ divisible. 
By the above, ${\cal J}/p^n{\cal J} \cong ({\cal I}_n)^m$ for a unique line 
bundle ${\cal I}_n$ and so the existence theorem applied to the 
${\cal I}_n$ show that there is an ${\cal I}$ with 
${\cal J} \cong {\cal I}^m$. 

Alternatively, the Kummer exact sequence shows that 
$\Pic(S)/m\Pic(S) \cong H^2_{et}(S,\mu_m)$ and 
$\Pic(\bar S)/m\Pic(\bar S) \cong H^2_{et}(\bar S,\mu_m)$ 
(here we use that $H^2_{et}(S,{\cal O}^*) = 0 = 
H^2_{et}(\bar S,{\cal O}^*)$) The result follows from proper 
base change (e.g. [Mi] p. 223).~\qed   

We need a variation of 1.2 where we have some control 
over values of functions at finitely many points. 
To this end, let 
$X$ be a scheme of finite type over 
a Noetherian ring $A$, and 
$P_1,\ldots,P_r$ a finite set of closed 
points each of which we write as $\iota_l: k(P_l) \to X$. In our 
application $X$ will be either $S$ or $\bar S$ 
so we will assume $X$ is projective over a Noetherian 
domain and is reduced. 
Form the sheaf ${\cal P}^* = \oplus_l \iota_l^*k(P_l)^*$. There is a 
surjective morphism of sheaves ${\cal O}^* = {\cal O}_X^* \to {\cal P}^*$ 
which is just evaluation and we let ${\cal O}^*_P$ 
be the kernel. Let ${\cal K}$ be the sheaf of total quotient 
rings of $X$ and ${\cal K}^*$ the group of units of ${\cal K}$. 
There are embeddings ${\cal O}_P^* \subset {\cal O}^* 
\subset {\cal K}^*$. Thus we have a exact sequence 
of sheaves $0 \to {\cal P}^* \to {\cal K}^*/{\cal O}^*_P \to 
{\cal K}^*/{\cal O}^* \to 0$. Since ${\cal P}^*$ and ${\cal K}^*$ 
are flasque 
we know $H^1(X,{\cal P}^*) = 0 = H^1(X,{\cal K}^*)$. 
Clearly $H^0(X,{\cal P}^*) = \oplus_l k(P_l)^*$. 
We set $K^* = H^0(X,{\cal K}^*)$. 
There are natural maps 
$K^* \to H^0(X,{\cal K}^*/{\cal O}^*_P)$ 
and $\oplus_l k(P_l)^* \to H^0(X,{\cal K}^*/{\cal O}^*_P)$, 
the later of which is an injection. 
The intersection, in $H^0(X,{\cal K}/{\cal O}^*_P)$, 
of the images of these maps we call $k^*$ and we can identify 
$k^*$ with the corresponding subgroup of $\oplus_l k(P_l)^*$. 
We have the exact diagram:
$$\diagram
\matrix{
&&0&&0\cr
&&\uparrow&&\uparrow\cr
K^*&\to&H^0(X,{\cal K}^*/{\cal O}^*)&\to&H^1(X,{\cal
O}^*)&\to&0\cr
||&&\uparrow&&\uparrow\cr
K^*&\to&H^0(X,{\cal K}^*/{\cal O}_P^*)&\to&H^1(X,{\cal
O}_P^*)&\to&0\cr
&&\uparrow&&\uparrow\cr
&&H^0(X,{\cal P}^*)&\to&H^0(X,{\cal P}^*)/k^*\cr
&&\uparrow&&\uparrow\cr
&&0&&0\cr}$$
Our goal is to interpret, a bit, $H^0(X,{\cal K}^*/{\cal O}^*_P)$ 
and $H^1(X,{\cal O}^*_P)$. Of course the former 
consists of equivalence classes of sets of pairs $\{(U_j,f_j )\}$ 
where $f_i \in {\cal K}^*(U_i)$, on $U_i \cap U_j$ 
the ratio $f_i/f_j$ is a unit, and this unit maps to 1 at all 
$P_l \in U_i \cap U_j$. 
If $\gamma = \{U_i,f_i\}$ is an element of $H^0(X,{\cal K}^*/{\cal O}^*)$ 
or $H^0(X,{\cal K}^*/{\cal O}^*_P)$ we say $\gamma$ 
{\bf avoids} ${\cal P}$ if for all $P_l$, all the relevant $f_i$ are 
units at $P_l$. 

Let $H^0_P(X,{\cal K}^*/{\cal O}^*_P)$, respectively 
$H^0_P(X, {\cal K}^*/{\cal O}^*)$ be the subgroup of those $\gamma$ 
which avoid all the $P_l$. The induced map 
$\rho: H^0(X,{\cal K}^*/{\cal O}^*_P) \to H^0(X,{\cal K}^*/{\cal O}^*)$ 
is onto and by definition $H^0_P(X,{\cal K}^*/{\cal O}^*_P)$ is the 
inverse image 
of $H^0_P(X,{\cal K}^*/{\cal O}^*)$. We need to prove 1.5 
but we begin with 1.4. 

\proclaim Proposition 1.4. Let $X$ be of finite 
type over a Noetherian ring $A$ with an ample bundle 
${\cal J}$. Fix an integer $m$ and a 
finite set of points $P_l$ on $X$. 
\smallskip 
a) Suppose $X = \P^r_A$. Then there is a homogeneous 
$f \in A[x_0,\ldots,x_r]$ of degree prime to $m$ and 
not $0$ at any $P_l$.  
\smallskip
b) There is a positive integer $r$ and a section $s$ of ${\cal J}^r$ 
such that $r$ is prime to $m$, ${\cal J}^r$ 
is very ample, and the support of $s$ contains none of the 
$P_l$. 
\smallskip
c) In particular, if $X$ is projective, there is an affine open $U \subset X$ 
containing all the $P_l$. 

\proof 
We begin with a). Let $Q_l \subset A[x_0,\ldots,x_r]$ 
be the homogeneous prime ideals associated to the $P_l$. 
Our argument will be the standard one, watching the degrees 
as we proceed. 
We induct on $s$, the cardinality of the set of $P_l$. 
If $s = 1$, we can take $f$ of degree 1. 
Assume the result for $s-1$. 
Choose $f_i$ of degree $d_i$, prime to $m$, such that $f_i \notin Q_j$ 
for $j \not= i$, $j = 1,\ldots, s$. We can assume $f_i \in Q_i$. 
Form $y = f_2^{t_2}{\ldots}f_s^{t_s}$ such that $d = d_2t_2 + \ldots + d_st_s$ is prime to $m$ and all $t_i > 0$. Note that 
$y \in Q_i$ for $i > 1$ and $y \notin Q_1$. 
Consider $f = f_1^d + y^{d_1}$, which has degree $d_1d$ prime to $m$. 
Then $f \notin Q_1$ because $y \notin Q_1$ and $f \notin Q_i$, $i > 2$ 
because $f_1 \notin Q_i$. Part a) is done. 

Next we claim there is an positive $r$, prime to $m$, 
such that ${\cal J}^r$ is very ample. This amounts to going 
through the proofs in [H] 7.6,  which uses arguments of 5.14 and 5.4 in [H],  
and being slightly careful. But now part b) reduces to a). 
Part c) is immediate.~\qed 

\proclaim Proposition 1.5. The maps $H^0_P(X,{\cal K}^*/{\cal O}^*) 
\to H^1(X,{\cal O}^*)$ and $H^0_P(X,{\cal K}^*/{\cal O}^*_P) 
\break \to H^1(X,{\cal O}^*_P)$ are surjective. 

\proof 
Suppose ${\cal I}$ is a divisor 
on $X$, viewed as an element of $H^0(X,{\cal K}^*/{\cal O}^*)$. 
We have assumed $X$ has an ample divisor ${\cal J}$. An easy argument from the 
definition (e.g. [H] p. 153) shows that ${\cal I} \otimes {\cal J}^n$ is ample 
for some $n$, and hence that ${\cal I}$ is the difference of ample divisors. 
Let ${\cal J}'$ be one of these ample divisors. 
By 1.4 there is a section of some ${\cal J'}^m$ 
whose support does not contain any of the $P_l$. 
Using 1.4 again, there is a section of ${\cal J'}^r$ 
whose support does not contain any of the $P_l$, where $r$ is prime to 
$m$. Using $a$ and $b$ such that $ar + bm = 1$, 
it is clear that each ${\cal J}'$ is represented by 
a class in $H^0(X,{\cal K}^*/{\cal O}^*)$ which 
misses all the $P_l$, and so the same applies to ${\cal I}$. 

That is, 
$H^0_P(X,{\cal K}^*/{\cal O}^*) \to H^1(X,{\cal O}^*)$ 
is surjective. Using the above diagram, it follows that 
$H^0_P(X,{\cal K}^*/{\cal O}^*_P) \to H^1(X,{\cal O}^*_P)$ 
is surjective.~\qed 

There is a well defined 
$\eta: H^0_P(X,{\cal K}^*/{\cal O}^*_P) \to \oplus_l k(P_l)^*$ 
given by evaluating the $f_i$ at the relevant $P_l$. 
It is immediate that $\eta$ is a splitting of the map 
$\rho: \oplus_l k(P_l)^* \to H^0_P(X,{\cal K}^*/{\cal O}^*_P)$ 
defined above. 
The inverse image of $H^0_P(X,{\cal K}^*/{\cal O}^*_P)$ 
in $H^0(X,{\cal K}^*)$ is $K^*_P$, defined as the subgroup 
of $K^*$ of all functions which are units at all the $P_l$. 
The following is now clear:

\proclaim Proposition 1.6. Let $K^*_P 
\!\subset\! H^0_P(X,{\cal K}^*/{\cal O}^*) \oplus [\oplus_l k(P_l)^*]$ 
via $g \to ((X,g), \sum_l g(P_l))$. Then $H^1(X,{\cal O}^*_P)$ 
is the quotient: 
$$H^0_P(X,{\cal K}^*/{\cal O}^*) \oplus [\oplus_l k(P_l)^*]\over K^*_P$$

Note that if 
$\gamma \in H^0_P(X,{\cal K}^*/{\cal O}^*)$ has support within 
a locally factorial open subset of $X$, then $\gamma$ 
can be identified with a (Weil) divisor whose support 
does not contain any of the $P_l$.  

We can now use 1.2 to show the following. Let $S \to \Spec(\Z_p)$ 
be as usual with $\bar S \subset S$ the reduced closed 
fiber. Assume $P_l$ are a finite set of closed points in 
$\bar S$ and $m$ is an integer prime to $p$.  

\proclaim Proposition 1.7. The canonical map induces an isomorphism 
 $${H^1(S,{\cal O}^*_P)\over m(H^1(S,{\cal O}^*_P))}\cong 
{H^1(\bar S,{\cal O}^*_P) \over m(H^1(\bar S,{\cal O}^*_P))}.$$

\proof 
Given the exact sequence $0 \to \oplus_l k(P_l)^*/k^* \to 
H^1(X,{\cal O}^*_P) \to H^1(X,{\cal O}^*) \to 0$ above, 
to prove this isomorphism it suffices to prove that 
$H^0(S,{\cal O}^*) \to H^0(\bar S, {\cal O}^*)$ 
is onto. But by [H] p. 277, 
$H^0(S,{\cal O}) \cong \lim H^0(S_n, {\cal O})$ 
where $S_n$ is the fiber of $p^n\Z_p$. 
Since units always lift modulo nilpotent ideals, 
we have the needed surjectivity.~\qed

\bigskip
\leftline{Section Two: Classification of Ramification}
\medskip 

In this section we assume $S$ is a nonsingular excellent surface. 
For any torsion abelian group $A$, $A' \subset A$ is the subgroup 
of elements of order prime to all residue characteristics of 
$S$. Let $K$ be the field of fractions of $S$ and $\alpha \in \Br(K)'$
an element of prime order $q$. We assume, for this section alone, 
that  $K$ contains a primitive $q$ root of one $\rho$, 
which we fix.
Using $\rho$, we define symbol classes etc. as in 
the introduction.  
For each curve $C \subset S$, the stalk ${\cal O}_{S,C}$ 
is a discrete valuation ring and so defines a ramification 
map $\Br(K)' \to H^1(F(C),\Q/\Z)' = \Hom_c(G_{F(C)},(\Q/\Z)')$ where 
$F(C)$ is the residue field of ${\cal O}_{S,C}$ and $G_{F(C)}$ is the Galois 
group of $F(C)$ in its separable closure. 
As in the introduction, elements of $\Hom(G_{F(C)},\Q/\Z)$ 
are identified with pairs $L/F(C),\sigma$ where $\sigma$ generates 
the Galois group of the cyclic extension $L/F(C)$. 
As observed above, the ramification locus of $\alpha$ 
is a finite union of curves on $S$. 
After blowing up (e.g. [L] p. 193), 
we can assume that that this ramification locus 
consists of nonsingular curves with normal crossings. 

What we study in this section includes the behavior of $\alpha$ 
with respect to {\bf all} the discrete valuation rings $R$ 
with $q(R) = K$ where $R$ lies over points or curves on $S$. 
Note that if $R$ lies over a curve of $S$, it equals ${\cal O}_{S,C}$ 
and so this is often not the hardest $R$ to understand. 
Thus let $P$ be a closed point of $S$, by which we mean a point 
of codimension 2. Let $R = {\cal O}_{S,P}$ be the stalk at $P$, 
which is a regular local ring of dimension 2. 
Let $M/K$ be a cyclic Galois extension of degree $q$. 
We will be most interested in results about when $M$ 
splits all the ramification of $\alpha$ over $R$.

We begin with a classification of the closed points of $S$ with respect 
to their relationship to the ramification locus of $\alpha$. 
Define $P \in S$ to be a {\bf distant} point if it is not on 
the ramification locus of $\alpha$. These points will rarely 
concern us. Define 
$P \in S$ to be a {\bf curve} point if it is on a single 
irreducible curve of the ramification locus. 
Finally, define $P \in S$ to be a {\bf nodal} point if it is 
a point in the intersection of two curves. 
It is the nodal points that will mostly require 
our analysis. If $u \in R'$, and $R'$ is a local ring, 
$\bar u$ is the image of $u$ in the residue field $R'/{\cal M}'$ of $R'$. 
Let us quote a result from [S] p. 32, slightly reworded and in 
our special case. 

\proclaim Theorem 2.1. Let 
$\alpha$ be as above, with ramification locus 
a union of nonsingular curves with normal crossings. 
If $C$ is a curve in that locus, let $L_C/F(C),\sigma_C$ 
be the ramification data of $\alpha$ at $C$. 
Let $R = {\cal O}_{S,P}$ be the stalk at a curve or nodal 
point $P$. 
In the following, $\alpha'$ 
always refers to an element of $\Br(R)$ and $u,v$ are always 
units in $R$. 
\smallskip 
a) If $P$ is a curve point and $C$ is the curve in the ramification 
locus containing $P$, then in $\Br(K)$, 
$\alpha = \alpha' + (u,\pi)$ where $\pi \in R$ is a prime defining 
$C$ at $P$. 
\smallskip
b) Suppose $P$ is a nodal point contained in both  
$C$ and $C'$ among the ramification locus of 
$\alpha$. Let $\pi$ and $\delta$ be primes of $R$ defining 
$C$, $C'$ respectively at $P$. Then either i) or ii) below hold: 
\smallskip
i) $\alpha = \alpha' + (u,\pi) + (v,\delta)$. 
\smallskip
ii) There is an $m$ prime to $q$ such that $\alpha = 
\alpha' + (u\delta^m, v\pi)$. 
\smallskip
Furthermore, the following holds. 
In a), $L_C/F(C)$ 
is unramified at $P$ and $\bar u^{1/p}$ defines $L_C/F(C),\sigma$ at that point. 
In b) i), $L_C/F(C)$ is unramified at $P$ and also defined by 
$\bar u^{1/q}$ at $P$. In b) ii), $L_C/F(C)$ is ramified at 
$P$ with ramification $m/q$ and defined by $(u\delta^m)^{1/q}$ 
at $P$. 

In all cases above, we call $\alpha - \alpha'$ a {\bf tail} 
of $\alpha$ at $R$ or $P$. 

We first consider the splitting at curve or distant points. 
The two cases are easy: 

\proclaim Theorem 2.2. If $P$ is a distant point, 
then $\alpha$ is unramified at any discrete valuation 
over $P$. Suppose $P$ is a curve point on $C$, 
and $C$ is in the ramification locus. Let 
$L/F(C), \sigma$ be the ramification data. 
If $L/F(C)$ splits at $P$, then $\alpha$ is unramified 
at any discrete valuation over $P$.  

\proof The distant point case is obvious. 
Let $P$ be a curve point on $C$. Write 
$\alpha = \alpha' + (u,\pi)$ where $u$ is a 
unit at $P$ with image $\bar u \in F(P)$. 
The residue field extension of $L/F(C)$ 
at $P$ is defined by $F(P)(\bar u^{1/q})$. 
That is, $L/F(C)$ splits at $P$ if and only if 
$\bar u \in (F(P)^*)^q$. Any valuation lying over 
$P$ will have $F(P)$ as a subfield of its residue field, and so 
it is obvious that $(u,\pi)$, and hence $\alpha$, 
is unramified at all such.~\qed  

It will be considerably more complicated to understand 
splitting all ramification at a curve point $P$ 
where $L/F(C)$ is not split. 
Let $C$ be a curve along which $\alpha$ ramifies and 
$P$ a nonsingular point on $C$. 
Let $R = {\cal O}_{S,P}$ and let $\pi = 0$ define $C$ at 
$P$. 
Write $\alpha = \alpha' + (u,\pi)$ as above. 
Set $F(P)$ to be the residue field of $R$. 
Suppose $L/F(C),\sigma$ is the ramification data of $\alpha$ at $C$. 
Suppose $x = \pi^s\delta \in R$ with $(s,q) = 1$ 
and $\delta$ is prime to $\pi$ in $R$. 
We are interested in when $M = K(x^{1/q})$ splits all the 
ramification of $\alpha$ over $R$.  
For convenience, may assume 
all the prime divisors of $\delta$ appear to prime to $q$ powers. 
To state the next result 
we successively blow up to form $\rho: S' \to S$ in such a way as 
to resolve 
the singularities in the (reduced) support of $x = 0$ at $P$. 
Let $\{E_i\}$ be the exceptional fibers of $\rho$. 
Write $(x) = \sum_i r_iE_i + \sum_j s_jC_j  + \sum_k t_kD_k$ 
where the $C_j$ are strict transforms of curves in $S$ containing $P$, 
and the $D_k$ are the curves in $S$ or $S'$ 
not containing $P$. We may take $C = C_1$ and (by definition) 
$s = s_1$. We call a curve or point {\bf relevant} if the 
residue field of that curve or point does not 
contain a $q$ root of $\bar u$. Of course we call a 
point or a curve {\bf irrelevant} if it is not 
relevant. 

\proclaim Theorem 2.3. Let $\alpha = \alpha' + (u,\pi)$ 
be as above, and assume $L/F(C)$ is nonsplit at a curve point $P$. 
Further assume we have blown up to resolve the singularities 
of $x = 0$ as above, and so the $r_i$ are defined.  
Then $M$ does not split the ramification 
of $\alpha$ at $P$ if and only if any of the $r_i$ 
for relevant $E_i$ are a multiple of $q$ 
or, baring that, any of the intersection points among the 
union of the $E_i$ and $C_j$ are relevant. 

\proof $L/F(C)$ is defined by $k(\bar u^{1/q})$ at $P$, 
and $\bar u$ is not a $q$ power in $F(P)$.
Thus $P$ itself is relevant. It follows that 
all the strict transforms $C_j$ are relevant. 
Also, by assumption, all the $s_j$ are prime to $q$. 
Suppose there is a relevant $E_i$ with $r_i$ a multiple 
of $q$. An irrelevant exceptional curve can only be created 
by blowing up an irrelevant point, and that is not $P$. 
That is, the first exceptional curve created is relevant. 
We can assume $E_i$ is the first relevant curve created 
in the resolution process with $r_i$ 
a $q$ multiple. 
Then $E_i$ arises from blowing up a relevant 
point on a relevant $E_{i'}$ with $r_{i'}$ prime to $q$. 
Thus at the end of the process $E_i$ will intersect $E_{i'}$ transversely 
at a relevant point $P'$ with $r_i$ and $r_{i'}$ as described 
and $P'$ being on no other curves in the support of 
$(x)$. Let $R_i = {\cal O}_{S',E_i}$. Then $M/K$ 
is unramified at $R_i$. Define 
$K_i/k(E_i)$ to be the residue field extension 
of $M/K$ at $E_i$. That is, $K_i = k(E_i)(\bar y^{1/q})$ 
where $\bar y$ is the image of some $y = x(z^q)$ 
and $y$ is a unit at $E_i$. It follows that 
$y$ has prime to $q$ valuation at $E_{i'}$, 
and hence that $K_i/k(E_i)$ ramifies 
at $P'$. Let $L_i/k(E_i),\sigma_i$ 
be the ramification data of $\alpha$ at $E_i$. 
Let $d_i$ be the discrete valuation corresponding 
to $E_i$. Since $d_i$ lies over $P$, 
the fact that $\alpha = \alpha' + (u,\pi)$ implies 
that $L_i = k(E_i)(\bar u^{1/q})$. Since $L_i$ 
is unramified at $P'$, it cannot equal $K_i$ 
and $M$ does not split the ramification of $\alpha$ 
with respect to $d_i$.  

Next assume all the relevant $E_i$ have $r_i$ prime to $q$ 
and $P'$ is a relevant intersection point. Then $P'$ is an intersection 
point of (say) local equations $\delta = 0$ and $\delta' = 0$ in the support 
of $(x)$. Both curves are relevant. If $R' = {\cal O}_{S',P'}$, then $x = 
w\delta^s\delta'^{t}$ where $w \in R^*$ and $s,t$ are prime 
to $q$. By 0.11 there is a valuation $d$ lying over $P'$ 
such that if $d(\delta) = a$ and $d(\delta') = b$ then 
$as + bt = nq$ and $q$ does not divide $ab$. Thus $M/K$ is unramified 
with respect to $d$. Since $P'$ lies over $P$, just as above 
the ramification 
of $\alpha$ at $d$ is $\bar u^{1/q}$. However, $M$ can be 
described as $K((w^{s'}\delta^b/\delta'^a)^{1/q})$ 
where $ss'$ is congruent to $b$ modulo $q$. By 0.11
it is clear that the residue field of $M$ does not contain 
$\bar u^{1/q}$ and once again we have a valuation where $M$ 
does not split the ramification of $\alpha$. 

Conversely, suppose all relevant $E_i$ have $r_i$ prime to 
$q$ and there are no relevant intersection points. 
Let $d$ be a valuation lying over the original $P$. 
Then $d$ must lie over a point or curve of the exceptional 
fiber. If $d$ lies over an irrelevant point or irrelevant 
curve, $\alpha$ is unramified at $d$. Thus we may assume $d$ 
lies over a relevant curve and since $M/K$ ramifies 
there, it follows that $M$ splits the ramification 
of $\alpha$ at any such $d$. Since $M/K$ is also ramified 
at ${\cal O}_{S,C}$, we are done.~\qed 

The main reason for stating and proving 2.3 was to show how complicated 
our analysis would have to be if we had to analyze extension 
fields $M = K(x^{1/q})$ that are as general as occur there. 
The following case is much simpler.

\proclaim Corollary 2.4. Suppose, in the situation 
of 2.3, $x = u\pi^s\delta^q$ in $R = {\cal O}_{S,P}$ 
where $u \in R^*$ and $s$ is prime to $q$. Then 
$M = K(x^{1/q})$ splits all the ramification of $\alpha$.  

\proof Here no blowing up is required and the result follows.~\qed  

We next classify what can happen at a nodal point $P$. 
Again set $R = {\cal O}_{S,P}$. 
If ${\cal M} \subset R$ is the maximal ideal, and $u \in R^*$, 
we let $\bar u$ be the image of $u$ in $F = R/{\cal M}$. 
If case b) ii) of 2.1 above holds we call $P$ a {\bf cold} point. 
We need to further 
analyze case b) i). Suppose $\bar u$, $\bar v$ do NOT generate the 
same subgroup of $F^*/(F^*)^q$. 
Then we say $P$ is a {\bf hot} point. If $\bar u,\bar v$ do generate 
the same subgroup of $F^*/(F^*)^q$, and they 
are not $q$ powers in $F$, we say $P$ is a 
{\bf chilly} point. 
If $1 \leq s \leq q-1$ is such that $\bar u^s\bar v^{-1} \in (F^*)^q$ 
we say that $s$ is the {\bf coefficient} of this chilly point with respect to 
$\pi$. Of course viewing the curves in the other order, 
if $s'$ is the coefficient of $P$ with respect to $\delta$, 
then $ss'$ is congruent to 1 modulo $q$.  
If both $u$, $v$ map to $q$ powers in $F$, we say $P$ is a 
{\bf cool} point.

The rest of this section will be a study of these four 
kinds of nodal points. We begin the the first of them. 

\proclaim Theorem 2.5. Suppose $P$ is a hot point. 
Then then the residual classes of $\alpha$ are not split 
by the ramification. In particular, $\alpha$ has index 
larger than $q$. 

\Demo Remark. The Jacob-Tignol example in [S] of an exponent 
$q$ and degree $q^2$ division algebra has a hot point, 
and the argument below is really theirs.  

\proof Write $\alpha = \alpha' + (u,\pi) + (v,\delta)$ 
as in 2.1 i). 
Since $\alpha$ ramifies at both $\pi$ and $\delta$, $u$ is not a 
$q$ power modulo $\pi$ and $v$ is not a $q$ power modulo $\delta$. 
We can assume the image of $u$ is not a $q$ power in $F = R/{\cal M}$. 
Let $R'$ be the localization of $R$ at $(\delta)$ with residue 
field $F'$  
and $L = K(\delta^{1/q})$. Let $\beta_{R'}$ 
be the residual Brauer class with respect to $L$.  
It is clear that $\beta_{R'} = \tilde \alpha' + (\tilde u,\tilde \pi)$ 
where the tilde refers to images in $\Br(F')$ and $F'^*$. 
Then $\tilde \pi$ defines a discrete valuation on $F'$, 
and with respect to this $\beta_{R'}$ has ramification $\bar u^{1/q}$, 
where $\bar u$ is the image of $\tilde u$ in $F$. 
The assumption that $P$ is a hot point implies 
that $F'(\tilde v^{1/q})$ does not split this ramification. 
But $\tilde v^{1/q}$ is the ramification of $\alpha$ at $\delta$, 
and we are done by 0.5.~\qed

Since in this paper we are concerned with division 
algebras of degree $q$, we often assume there are no 
hot points. Our next observation is that we can blow up 
to eliminate any cool points. 

\proclaim Theorem 2.6. Suppose $P \in S$ is a cool point. 
Then if we blow up $S$ at $P$, the Brauer group 
element $\alpha$ does not ramify on the exceptional divisor, 
and so the cool point has been turned into two curve points. 

\proof Let $R,{\cal M}$ be the local ring of $S$ at $P$, a cool point. 
Then a tail of $\alpha$ can be chosen to look 
like $[(u,\pi)_q] + [(v,\delta)_q]$ where $\bar u,\bar v$ 
are $q$ powers in $F = R/{\cal M}$. If $R'$ 
is a discrete valuation lying over ${\cal M}$ with valuation $d$, 
then the residue of this tail has the form 
$\bar u^{d(\pi)}\bar v^{d(\delta)}$ and so 
is a $q$ power in $F$, which is a subfield 
of the residue field of $R'$. That is, 
$\alpha$ is unramified at every discrete valuation 
over ${\cal M}$, implying it is unramified on the exceptional 
divisor.~\qed 

For the rest of this section we will assume we have used 2.6 to eliminate 
any cool points and that there are no 
hot points. Note that this means the following. 
Let $P$ be an intersection point 
of two curves $C$, $C'$ along which $\alpha$ 
ramifies with covers $L/F(C)$ and $L'/F(C')$. 
Then either $P$ is a ramified point with respect to both 
extensions or $P$ is a nonsplit point with respect to both 
extensions. 
We are left with studying chilly and cold points. 
Let us begin with chilly points. 

\proclaim Proposition 2.7. Suppose $P$ is a chilly point, 
$R = {\cal O}_{S,P}$  
and $\pi \in R$, $\delta \in R$ are the two primes defining 
the ramification locus of $\alpha$ at $P$. Let 
$s$ be the coefficient with respect to $\pi$, and $w$ a 
unit of $R$. 
\smallskip 
a) $M = K((w\pi\delta^s)^{1/q})$ 
splits all the ramification of $\alpha$ at any prime lying over 
$R$. 
\smallskip
b) For any $t$ not congruent to $s$ modulo $q$, 
$M' = K((w\pi\delta^t)^{1/q})$ fails to split the ramification 
of $\alpha$ at some prime lying over $R$. 

\proof Suppose $M$ is as described in a). Let 
$d: M \to \Z$ be a valuation lying over $R$. 
If $d$ lies over any height one prime not $\pi$ or $\delta$, 
or if $d(\pi)$ and $d(\delta)$ are both $q$ multiples, 
then clearly $\alpha$ is unramified at $d$. If $d$ lies over 
$\pi$ or $\delta$, then $M/K$ is ramified at $d$ and so  
$M$ splits the ramification of $\alpha$ at $d$. 
Thus we may assume $d$ lies over the 
maximal ideal, ${\cal M}$, of $R$, $d(\pi) = a > 0$ and $d(\delta) = b > 0$, 
and one of the $a,b$ is prime to $q$. 
Note this also means that $k = R/{\cal M}$ is a subfield of the residue field 
of $d$. The ramification of 
$(u,\pi)$ at $d$ is $\bar u^{a/q}$ and the ramification of 
$(v,\delta)$ is $\bar v^{b/q} = \bar u^{bs/q}$, 
so the ramification of $\alpha$ is $\bar u^{(a + bs)/q}$. 
If $a + bs$ is prime to $q$, then $M/K$ is ramified at $d$ 
and so splits the ramification of $\alpha$. If $a + bs$ 
is a multiple of $q$, $\alpha$ is not ramified at $d$ and we are done. 

Continuing with b), let $M'$ be as defined and $k = R/{\cal M}$. 
By 0.11 there is 
a valuation $d$ on $K(T)$ lying over $R$ with 
the following properties. First of all, 
$d(T) = 1$, and $d(\pi) + td(\delta) = mq$.  
Secondly, the residue field of $d$ is $k(\pi',\delta')$  
where $\pi' = \pi/T^{d(\pi)}$, $\delta' = \delta/T^{d(\delta)}$. 
Note that $x = w\pi\delta^t/T^{mq} = \pi'\delta'^t$ has image 
$\bar x$ which is  
part of a transcendence base $\bar x, z$ of $k(\pi', \delta')$ over $k$. 
Since $\bar u^t\bar v$ is not a $q$ power in 
$k$, and $M'(T)$ has residue field $k(x^{1/q}, z)$ with respect to the 
unique extension $d''$, of $d$, it follows that the ramification of $\alpha$ 
is not split at $d''$ in $M'(T)$, and hence not split by the 
restriction of $d''$ to $M'$.~\qed 

Besides the splitting question handled above, we will need 
some results about the residual Brauer class in case 
a) above. 

\proclaim Theorem 2.8. Suppose $P$ is a chilly point 
at the intersection of $C$ and $C'$ in the ramification 
locus of $\alpha$. Let $C$, $C'$ be 
locally defined by $\pi = 0$ and $\delta = 0$ 
respectively and let $s$ be the  
coefficient with respect to $C$. 
Let $M = K((w\pi\delta^s)^{1/q})$ as in 2.7 a) 
above. Suppose $\beta_C$ and $\beta_{C'}$ 
are the residual Brauer classes of $\alpha$ with respect to $M/K$. 
Then $\beta_C$ and $\beta_{C'}$ are both unramified 
at $P$ and have equal images in $\Br(F(P))$. 

\proof Let $ss' - 1$ be divisible 
by $q$, so $s'$ is the coefficient with respect 
to $C'$. We can also write $M = 
K((w^{s'}\delta\pi^{s'})^{1/q})$. 
At $R = {\cal O}_{S,P}$ write $\alpha = 
\alpha' + (u,\pi) + (v,\delta)$ where $\alpha' 
\in \Br(R)$, $u,v \in R^*$, and $u^s$ and $v$ differ by 
$q$ powers in $F(P)$. Denote by $L/F(C),\sigma$ 
and $L'/F(C'),\sigma'$ the ramification 
data of $\alpha$ at $C$ and $C'$. Then $L/F(C)$ is defined 
by $\bar u^{1/q}$ at $P$ and $L'/F(C')$ is defined by 
$\bar v^{1/q}$ at $P$. The image of $\alpha$ in $\Br(M)$ 
is the same as the image of $\alpha'' = \alpha' + (u,w^{-1}\delta^{-s}) + 
(v,\delta) = \alpha' + (u,w^{-1}) + (v/u^s,\delta)$. 
Since $v/u^s$ is a $q$ power at $P$, the image of 
$\alpha''$  in $\Br(F(C))$ is unramified at $P$. 
Moreover the image of $\alpha''$ in $\Br(F(P))$ 
is $\bar \alpha' + (\bar u,\bar w^{-1})$. 
Looking at $\beta_{C'}$, which means 
reversing $\pi$ and $\delta$, and therefore 
switching $s$ and $s'$ and $u,v$, we get the image 
$\bar \alpha' + (\bar v,\bar w^{-s'})$ which is the same.~\qed 

Ultimately, we are going to show $\alpha$ is cyclic by 
finding an element $f$ where the support of $(f)$ includes the 
full ramification locus of $\alpha$ and the coefficients of 
$(f)$ are chosen so that a) above applies and not b). There is an 
inherent difficulty with this if there are ``loops'' of curves 
where incompatible coefficients are required to meet 
condition a) above. To get around this, we consider the 
effect of blowing up on a chilly point. 

Let $[(u,\pi)_q] + [(v,\delta)_q]$ be a tail 
of $\alpha$ at $R = {\cal O}_{S,P}$ with coefficient 
$s$ with respect to $\pi$. The blowup defines a valuation 
with $d(\pi) = d(\delta) = 1$, and so the ramification of $\alpha$ at the 
blowup is $\bar u\bar v$ which is the same as $\bar v^{s+1}$ 
modulo $q$ powers. Thus if $s + 1$ is a multiple of $q$, 
there is no ramification on the blowup and we have turned 
a chilly point into two curve points. In any other case, 
there are 2 nodal points to consider. Let $R'$ be the local 
ring at the intersection of the strict transform $\pi = 0$ and the blowup. 
Then in $R'$ we have a $\zeta$ with $\zeta\delta = \pi$ 
where $\zeta = 0$ defines the strict transform of $\pi = 0$ 
and $\delta = 0$ defines the blowup divisor. 
Thus the tail of $\alpha$ at $R'$ is $(u,\zeta) + (uv,\delta)$. 
It follows that $R'$ is a chilly point with coefficient $s+1$ 
with respect to $\zeta = 0$. 
Similarly, let $R''$ be the intersection of the blowup  
with the strict transform of $\delta = 0$ and let $s'$ 
be the coefficient of $P$ with respect to $\delta$. The same argument 
shows that if $P''$ is the intersection of the blowup with 
$\delta = 0$, the coefficient is $s' + 1$ at that point.

Consider a graph whose vertices are the curves 
in the ramification locus, and the edges are the chilly points. 
Two vertices have an edge between them if they both 
contain that chilly point. For any edge, 
blowing up can have one of two effects. If the coefficient 
is $q - 1$, blowing up removes the edge. Otherwise, blowing 
up adds a vertex between the two vertices and two edges 
connecting the new vertex with both of the old ones. 
A loop in the above graph we call a {\bf chilly loop}. 
It is clear that by repeated blowing up we can break any chilly loop. 

\proclaim Corollary 2.9. After repeated blowing up, we can assume 
there are no chilly loops in the ramification locus of 
$\alpha$. 

\proclaim Corollary 2.10. Suppose $C_i$ are all the 
curves in the ramification locus and we have blown up 
so that there are no chilly loops. 
Then we can choose, for each $C_i$, a nonzero
$s_i \in \Z/q\Z$ such that the following holds.  
Suppose $P$ is a chilly point on $C_i$ and $C_j$ with coefficient $s$ 
with respect to $C_i$. Then 
$s = s_j/s_i$ in $\Z/q\Z$. 

\proof The graph is a tree so this is an easy induction,  
one leaf at a time.~\qed 

It now behooves us to consider splitting 
at cold points. More specifically, suppose 
$P$ is a cold point defined locally by the intersection 
of curves $C$ and $C'$ in the ramification locus of $\alpha$. 
Let $R = {\cal O}_{S,P}$ and let $\pi$ and $\delta$ be primes 
of $R$ defining $C$ respectively $C'$ at $P$. 
Suppose $s,t$ are prime to $q$. We are interested in when 
$M = K((w\pi^s\delta^t)^{1/q})$ splits all the ramification 
of $\alpha$ over $R$. What we will find is that this is 
determined by the residual Brauer class $\beta_C$ of 0.5. 
Recall that $\beta_C \in \Br(F(C))$, where $F(C)$ 
is the residue field of $R$ localized at $C$. 
By assumption,  
$P$ is nonsingular on $C$ so defines a discrete valuation 
on $F(C)$. That is, if $F(P)$ is the residue field at $P$, 
$\beta_C$ has some ramification $\chi_P \in \Hom(G_{F(P)},\Q/\Z)$. 

Our immediate goal is a second description of $\chi_P$ in terms 
of ramification on $K = F(S)$. Let $d'$ be a discrete valuation 
of $K$ lying over $P$, and set $a = d'(\pi)$ and $b = d'(\delta)$. 
Let $s'$ be the inverse of $s$ modulo $q$. 
Assume $M/K$ is unramified at $d'$, which is equivalent to 
assuming $sa + tb$ is divisible by $q$. Let $d$ be any extension 
of $d'$ in $M$. Note that $F(P)$ is a subfield of the residue 
field of $M$ at $d$.  

\proclaim Proposition 2.11. Suppose $P$ is a cold 
point and $M = K((w\pi^s\delta^t)^{1/q})$, 
$\beta_C$, $\chi_P$, $d$ are as above. The ramification of $\alpha$ 
at $d$ is the image of $\chi_P^{b}$. 

\proof 
By 2.1 we can write $\alpha = \alpha' + (u\delta^m,v\pi)$ 
for $m$ prime to $q$. Then 
$\alpha$ has the same image in $\Br(M)$ as $\alpha'' = 
\alpha' + (u^m\delta^m,vw^{-s'}\delta^{-s't})$ which is manifestly 
unramified with respect to $C$ and so has image $\beta_C$ 
in the residue field. In addition, $\alpha'$ maps to an element 
of $\Br(F(C))$ unramified at $P$. Finally, the image, $\bar \delta$, of 
$\delta$ in $F(C)$ is the prime defining $P$, and the images, 
$\bar u$, $\bar v$, $\bar w$,  of $u$, $v$, $w$ are all units 
at $P$. All together, $\chi_P$ is defined by $x^{1/q}$ 
where $x$ is the image of $(\bar u^{m(-s't)}/\bar v^m(\bar w^{-s'm}))$ 
which up to $q$ powers is $(\bar w/(\bar u^t\bar v^s))^{s'm}$. 
On the other hand, the ramification of $\alpha$ with respect to 
$d$ is the image of the ramification of $\alpha''$ 
with respect to $d'$ and this (by the formula) is $y^{1/q}$ 
where $y$ is the image of $(w/(u^tv^s))^{bs'm}$.~\qed 

We can use the above calculations to observe 
a relationship between the ramification of the residual 
classes at cold points.  

\proclaim Corollary 2.12. Suppose 
$P$ is a cold point at the intersection 
of $C$, $C'$ in the ramification locus. Let $M = 
K((w\pi^s\delta^t)^{1/q})$ be as above.  
Then the ramification of $s\beta_C$ and $-t\beta_{C'}$ 
are equal at $P$. 

\proof Of course we have fixed a $q$ root of one $\rho$, 
and it is easy to see that our description of the 
tail of $\alpha$ implies that $L/F(C),\sigma$ has ramification 
$\rho^{m'}$ at $P$, where $mm' - 1$ is divisible by $q$. 
By the proof of 2.11, and using $\rho$ again, 
the ramification of $\beta_C$ is represented by $x^{1/q}$ 
where $x$ is the image of $(w/u^tv^s)^{ms'}$ and where $ss' - 1$ 
is divisible by $q$. 
To reverse the roles of $C$ and $C'$ we can also write 
$\alpha = \alpha' + (v^{-m}\pi^{-m},u\delta)$.  
If $L'/F(C')$ is the ramification of $\alpha$ 
at $C'$, then $L'/F(C')$ has ramification 
$\rho^{-m'}$ at $P$. The same argument as in 2.11 
shows that $\beta_{C'}$ has ramification $x'^{1/q}$ where 
$x'$ is the image of $(w/u^tv^s)^{-mt'}$.~\qed 

The ramification of $\beta_C$ at a cold point determines 
the splitting of the ramification of $\alpha$ at that point:

\proclaim Corollary 2.13. Suppose $P$ is a cold point defined 
locally as the intersection of $C$ and $C'$ 
in the ramification locus. Let   
$R = {\cal O}_{S,P}$ and $M = K((w\pi^s\delta^t)^{1/q})$, for some
$s$, $t$ prime to $q$. Let $\beta_C$ be the residual Brauer class 
of $\alpha$ with respect to $M$. Then $M$ 
splits all the ramification of $\alpha$ over $R$ if and 
only if $\beta_C$ is unramified at $P$.  

\proof 
Let $d: K^* \to \Z$ be a valuation over $R$ at which $\alpha$ 
ramifies. If $d$ lies over a prime of height one, 
it must be $\pi$ or $\delta$ and $M$ is ramified at those 
primes. Thus we may assume $d$ lies over the maximal ideal 
of $R$. Let $d(\pi) = a > 0$ and $d(\delta) = b > 0$. 
If both $a,b$ are divisible by $q$, $\alpha$ does not 
ramify at $d$, so we assume one of $a$ or $b$ is prime 
to $q$. If $sa + tb$ is not divisible by $q$, then 
$M/K$ ramifies at $d$. Thus we may assume $M/K$ 
is unramified at $d$. If $\chi_p$ 
is trivial, 2.11 shows that $M$ splits the ramification 
at any such $d$.  

Conversely, by 0.12, there is a valuation $d$ on $K(T)$ 
where $sd(\pi) + td(\delta)$ is divisible by $q$ and the residue field 
of $d$ is $F(P)(\pi',\delta')$ as described there. 
If $M$ splits the ramification 
of $\alpha$ at the restriction of that $d$, then $M(T)$ 
must split the ramification of $\alpha$ at $d$. 
In the notation of 2.11 it follows that $(w/u^tv^s)$ must 
map to a $q$ power in $F(P)(\pi',\delta')$ from which 
the result is clear.~\qed 

While on the subject of residual Brauer classes, for completeness 
we add: 

\proclaim Corollary 2.14. Suppose $P$ is a curve point on $C$ 
and the ramification $L/F(C)$ splits at $P$. 
Suppose $M = K(\pi^{1/q})$ and $\pi$ has $C$ valuation 
prime to $q$. If $\beta_C$ is the residual Brauer class with respect 
to $C$, then $\beta_C$ is unramified at $P$. 

\proof Let $R = {\cal O}_{S,P}$. 
We can write $\alpha = \alpha' + (u,\pi_C)$ where $\pi_C = 0$ 
defines $C$ locally at $P$. Also, $\pi = v\pi_C^s\delta$ where 
$v \in R^*$, $s$ is prime to $q$, and $\delta$ is not divisible 
by $\pi_C$. $\alpha$ has the same image in $\Br(M)$ as 
$\alpha'' = \alpha' + (u,v''\delta')$ where $v'' \in R^*$ and $\delta''$ 
is not divisible by $\pi_C$. Since $L/F(C)$ is split at $P$, 
the image, $\bar u$, of $u$ in $F(P)^*$ is a $q$ power. 
Direct calculation shows that since $\beta_C$ is the image of $\alpha''$, 
$\beta_C$ is unramified at $P$.~\qed 

\bigskip
 
\leftline{Section Three: Adding a root of unity}
\medskip 
The purpose of this section is to detail how the results 
of section two can be extended to the case where 
$K = F(S)$ does not contain a primitive $q$ root of one. 
To this end, we begin more generally. 

Let $R$ be a regular local ring of dimension 2 with 
residue characteristic $p$, maximal ideal ${\cal M}$ 
and fraction field $K$. 
Let $m$ be an integer prime to 
$p$. Let $\mu_m$ be the group of $m$ 
roots of one over $K$ with generator $\rho \in \mu_m$.  
Let $f(x)$ be the monic minimal polynomial of $\rho$ 
over $K$, so $f(x) \in R[x]$ and we can set $R' = R[x]/(f(x))$. 
Then $R'/R$ is Galois with (abelian) group $H$. 
For any prime of $R$, the group $H$ acts transitively on the 
primes of $R'$ lying over $R$. 

Assume $\pi,\delta \in {\cal M} - {\cal M}^2$ and ${\cal M} = (\pi,\delta)$. 
Let $H_{\cal M}$, $H_{\pi}$, and $H_{\delta}$ 
be the stabilizers of one (and hence all) of the 
prime ideals lying over  ${\cal M}$, $(\pi)$, or $(\delta)$ 
respectively. If $R/{\cal M}$ contains a primitive $m$ 
root of one, $H_{\cal M} = 1$. By 0.12 $R'$ is a UFD:  

\proclaim Lemma 3.1. One can choose $\pi_1$ such that $\pi_1$ generates a 
prime over $(\pi)$ and such that the stabilizer 
of $\pi_1$ as an element is $H_{\pi}$.  
A similar result (of course) holds for $(\delta)$. 

\proof 
$R'^{H_{\pi}}$ is a UFD by 0.12. 
Let $\pi_1$ generate a prime of $R'^{H_{\pi}}$ lying over $\pi$.~\qed 

By 3.1, we can write 
$\pi = u\prod_i \pi_i$ and $\delta = v\prod_j \delta_j$ 
to be the prime decompositions of $\pi$ and $\delta$ in 
$R'$ (where $u,v \in R'^*$). It is immediate that all the 
$(\pi_i)$ and $(\delta_j)$ are distinct and each set forms a
single $H$ orbit. By 3.1 b), we can assume the sets of 
elements $\{\pi_i\}$ and $\{\delta_j\}$ form a single 
$H$ orbit. Hence by changing our choice of $\pi$, $\delta$ 
we can assume $\pi = \prod_i \pi_i$ and $\delta = \prod_j \delta_j$. 
This is merely a convenience. 

Let $\{{\cal M}_k\}$ be the set of maximal ideals of $R'$ and 
${\cal J} = \cap _k {\cal M}_k$ the Jacobson radical. Since $R'/R$ 
is etale, ${\cal J} = (\pi,\delta)R'$. Any ${\cal M}_k$ contains 
$\pi$ and $\delta$ and hence at least one $\pi_i$ 
and one $\delta_j$. 
Since ${\cal J}R'_{{\cal M}_k} = {\cal M}_kR'_{{\cal M}_k}$ we have 
$(\pi,\delta)R'_{{\cal M}_k} = (\pi_i,\delta_j)R'_{{\cal M}_k} = 
{\cal M}_kR'_{{\cal M}_k}$. 
If $\pi_i$ and $\pi_{i'}$ were in the same ${\cal M}_k$, 
then $\pi \in {\cal M}_k^2$ which would contradict the above. 
Thus each ${\cal M}_k$ contains a unique $\pi_i$ and $\delta_j$. 
However, multiple ${\cal M}_k$ can contain the same $\pi_i$ and $\delta_j$. 
Checking locally, it follows that $(\pi_i,\delta_j)$ 
is the intersection of a uniquely defined set of maximal ideals 
of $R'$.  

\proclaim Lemma 3.2. Let $R'/R$, $\pi = \prod_i \pi_i$ and 
$\delta = \prod_j \delta_j$ be as above. 
For a fixed ${\cal M}_k$, there is a unique $\pi_i$ and a unique 
$\delta_j$  in ${\cal M}_k$. 
In particular, $H_{\cal M} \subset H_{\pi} \cap H_{\delta}$. 
The ideal $(\pi_i,\delta_j)R'$ is either $R'$ or is the 
intersection of maximal ideals of $R'$ which form a single and 
unique $H_{\pi} \cap H_{\delta}$ orbit.

\proof The inclusion $H_{\cal M} \subset H_{\pi} \cap H_{\delta}$ 
is immediate from the uniqueness of $\pi_i$, $\delta_j$ 
in some ${\cal M}_k$. Looking locally, it is clear that 
$(\pi_i,\delta_j)$ is the intersection of the maximal ideals 
containing it. This set of maximal ideals is clearly 
closed under the action of $H_{\pi} \cap H_{\delta}$. 
Assume $(\pi_i,\delta_j) \not= R'$. 
If $R'' = R'^{H_{\pi} \cap H_{\delta}}$, then 
$R''/R$ is Galois with group $\bar H = H/(H_{\pi} \cap H_{\delta})$. 
Every maximal ideal of $R''$ has trivial stabilizer. 
If $h \in H$ is not in $H_{\pi} \cap H_{\delta}$, then 
$h(\pi_i) \not= \pi_i$ or $h(\delta_j) \not= \delta_j$. 
Thus $h(\pi_i,\delta_j)R' \not= (\pi_i,\delta_j)R'$. 
It follows from a counting argument 
that $(\pi_i,\delta_j)R''$ is contained in 
a unique maximal ideal of $R''$, and so $(\pi_i,\delta_j)R''$ is a 
maximal ideal of $R''$. The rest of the lemma is now 
immediate.~\qed 

Let $S$ be a nonsingular excellent surface. Let $S' \to S$ be the Galois 
cover gotten by adjoining a primitive $q$ root of one, 
and let $K' = F(S')$ be the function field of $S'$. 
We assume this extension is etale, meaning that $q$ is prime to 
all the residue characteristics of $S$. 
Let $H$ be the Galois group of $S'/S$, so $H$ is cyclic of order 
$m$ dividing $q-1$. The fact that $m$ 
is prime to $q$ is behind much of this section. 
For any curve $C \subset S$, or point 
$P \in S$, let $H_C$ or $H_P$ be the stabilizer of one and hence 
any point or curve lying over $P$ and $C$ respectively. 

We are interested in applying section two to $S'$ 
as a way of classifying ramification on $S$. 
To this end, we rephrase 3.1 and 3.2 above. 
Let $P$ be a point on $S$ 
which is the local scheme theoretic intersection of nonsingular curves 
$C$ and $C'$, said curves locally defined by $\pi = 0$ 
and $\delta = 0$. Then by 3.1 and 3.2, the points of $S'$ 
mapping to $P$ are each locally scheme theoretically defined 
by a unique $C_i$ and $C_j'$ in $S'$, where the $C_i$ lie over 
$C$ and the $C_j'$ lie over $C'$. In addition, none of the 
$C_i$ intersect each other in $S'$, and similarly for the $C_j'$. 

Fix an element $\alpha \in \Br(F(S))$ of order $q$. Then, 
as usual, $\alpha$ has a ramification locus which is a 
bunch of curves $C_i$ and cyclic covers $L_i/F(C_i), \sigma_i$ 
where $\sigma_i$ generates the Galois group of $L_i/F(C_i)$.  
Also as usual, we can blow up and assume the $C_i$ are all nonsingular 
with normal crossings. It will be important for us to understand the 
relationships 
between this ramification data and the corresponding 
data over $S'$. To this end, let $\alpha'$ be the image of 
$\alpha$ in $\Br(F(S'))$ and $L_j'/F(C_j'),\sigma_j' $ 
the ramification data of $\alpha'$. 

\proclaim Theorem 3.3. The $C_j'$ are precisely the preimages of 
the $C_i$ in $S'$. If $C_j'$ lies over $C_i$, then $L_j' = 
L_i \otimes_{F(C_i)} F(C_j')$. $\sigma_j'$ is the extension of 
$\sigma_i$ trivial on $F(C_j')$. Furthermore, for 
fixed $i$, none of the inverse images of $C_i$ 
intersect.  

\proof This is obvious from 0.3, noticing that 
$S'/S$ is unramified everywhere.~\qed 

We will parallel section two  
and classify the points of $S$ with respect to the ramification 
data of $\alpha$. 
The easy things are easy. If $P$ is not on any $C_i$, we say 
$P$ is a {\bf distant} point. If $P$ is on exactly one $C_i$, 
we say $P$ is a {\bf curve} point. It is obvious from 3.3 that $P$ 
is a distant or curve point if and only if one and hence 
all of its preimages in $S'$ have the same behavior with respect to 
$\alpha'$. It is also obvious that $\alpha$ is unramified at any 
discrete valuation over a distant point. 
The following is obvious from 2.2.  

\proclaim Theorem 3.4. Suppose $P$ is a distant point, 
or a curve point where the ramification $L/F(C)$ 
is split. Then $\alpha$ is unramified at any discrete valuation 
over $P$. 

We also need to generalize (trivially) 2.4.  

\proclaim Proposition 3.5. Suppose $P$ is a curve point 
on $C$, $R = {\cal O}_{S,P}$, and $x = u\pi^s\delta^q$ 
where $u \in R^*$ and $s$ is prime to $q$. 
Then $M = K(x^{1/q})$ splits all the ramification of 
$\alpha$ over $P$. 

If $P$ is on exactly two of the $C_i$, we say $P$ is a 
{\bf nodal} point. Then all the preimages of 
$P$ are nodal points. Suppose $P$ is on $C_1$ and $C_2$ 
and $L_k/F(C_k),\sigma_k$ is the ramification of $\alpha$ 
at $C_k$ for $k = 1,2$. Assume $P'$ is a point on $S'$ mapping 
to $P$ and $C_1'$, $C_2'$ are curves on $S'$ which lie over 
$C_1$, $C_2$ and both contain $P'$. Then $L_k/F(C_k)$ 
is ramified at $P$ if and only if $L_k' = L_k \otimes_{F(C_k)} F(C_k')$ 
is ramified at $P'$. If $L_k/F(C_k)$ is unramified at $P$, 
then $P$ splits in $L_k$ if and only if $P'$ splits in $L_k$. 
Finally, if $P$ extends uniquely in $L_k$, then the 
residue field of $L_k'/F(C_k')$ at $P'$ is the extension 
of that of $L_k/F(C_k)$ at $P$. 

\proclaim Theorem 3.6. Let $P$ be a nodal point. Then  every 
preimage of $P$ in $S'$ is a nodal point for $\alpha'$. 
If one of the preimage points of $P$  is hot, or chilly, or cool, 
or cold, then all have the identical behavior. 
Furthermore, if all the preimages of $P$ are chilly, then they all 
have the same coefficient with respect to the corresponding preimage 
of $C_1$. 

\proof 
The result for cold points is obvious. In all other cases, the 
definitions of section two looked at elements $\bar u, \bar v \in F(P)^*$. 
It suffices to use the observation of 2.1 
that $F(P)(\bar u^{1/q})$ and $F(P)(\bar v^{1/q})$ 
are the residue extensions at $P$ of the respective 
ramification extensions $L_1/F(C_1)$ and $L_2/F(C_2)$.~\qed 

\Demo Definitions. Clearly, then, it makes sense 
to say a point $P$ of $S$ is {\bf hot}, or {\bf chilly}, or 
{\bf cool}, or {\bf cold} 
if one and hence all its preimages have the same property. 

It will also be useful to rephrase the condition of being a chilly  
point with coefficient $s$ with respect to 
$C_1$. Suppose $L_1/F(C_1),\sigma_1$ and $L_2/F(C_2),\sigma_2$ 
is the ramification data for $\alpha$ at a nodal point 
$P$.  

\proclaim Corollary 3.7.  
a) Suppose $P$ is a chilly point with coefficient $s$ with respect to 
$C_1$. 
Then both $L_i$ are unramified at $P$. If $\bar L_i/F(P)$ 
are the induced residue extensions, then both are fields  
unramified at $P$ and are equal. If $\bar \sigma_1$ and $\bar \sigma_2$ 
are the induced generators of the Galois groups, then 
$\bar \sigma_2^s = \bar \sigma_1$. 
\smallskip
b) Suppose $P$ is a cold point. Then both $L_1/F(C_1),\sigma_1$ 
and $L_2/F(C_2), \sigma_2$ are ramified at $P$. If $\rho$ 
is the ramification of $L_1/F(C_1),\sigma_1$ at $P$, 
then $\rho^{-1}$ is the ramification of $L_2/F(C_2),\sigma_2$ at $P$.

\proof 
Both parts are proven by extending to $K'$ and using the fact 
$K'/K$ has degree prime to $q$. In a), we have just rewritten 
the definition of chilly and coefficient. In b), 
we note the same fact at cold points in $K'$ and again translate.~\qed

There are a further series of consequences of 3.6. 

\proclaim Proposition 3.8. If we blow up a cool 
point, then this point is replaced in the ramification 
locus by two curve points. After repeated blowing up, 
we can assume there are no chilly loops. With this, if $C_i$ 
form the ramification locus of $\alpha$ then for each 
$C_i$ we can choose a nonzero $s_i \in \Z/q\Z$ 
such that the following holds. Suppose $P$ is a chilly point 
which is the locally the intersection of $C_i$ and $C_j$, 
and which has coefficient $s$ with respect to $C_i$. 
Then $s = s_j/s_i$ in $\Z/q\Z$. 

\proof The blow up of a point $P$ on $S$ 
pulls back to the successive blow up (any order) of the 
preimage points. This makes the rest of the proposition 
clear. Since there are no chilly loops, the last sentence is 
clear just as in 2.10.~\qed  

Next we turn to generalizing 2.7. That is, we consider the 
coefficient at a chilly point and the consequences for splitting. 

\proclaim Proposition 3.9. Let $P$ be a chilly point and 
$\pi = 0$ and $\delta = 0$ the local equations for the two curves 
through $P$ along which $\alpha$ ramifies. Let $s$ be the 
coefficient with respect to $\pi$. 
\smallskip 
a) $L = K((\pi\delta^s)^{1/q})$ 
splits all the ramification of $\alpha$ at any prime lying over 
$R$. 
\smallskip
b) For any $t$ not congruent to $s$ modulo $q$, 
$L_t = K((\pi\delta^t)^{1/q})$ fails to split the ramification 
of $\alpha$ at some prime lying over $R$. 

\proof Obviously we will extend scalars to $K' = K(\mu_q)$ 
and use 2.7 to prove this. Let $R = {\cal O}_{S,P}$ and $R'/R$ 
the extension gotten by adjoining $\mu_q$. Let $\pi = \prod_i \pi_i$ 
and $\delta = \prod_j \delta_j$ be the prime decompositions 
in $R'$. At each closed point of $R'$ defined by $(\pi_i, \delta_j)$, 
$L$ and $L_t$ have the form $K((wu_i\delta_j^s)^{1/q})$ and 
$K((w\pi_i\delta_j^t)^{1/q})$ respectively for a unit $w$ 
at that point. Thus 2.7 applies. In a), we conclude that 
$L' = L \otimes_K K'$ splits all the ramification over all points 
over $P$. Since $L'/L$ has degree prime to $q$, we are done in a). 
In b), we find the discrete valuation over some preimage point 
that fails to kill the ramification and restrict it to $K$.~\qed 

We need to make some remarks about how adding roots of one 
effects residual Brauer classes. 
Suppose $P \in S$ is a nodal point and the intersection 
of $C$ and $C'$ in the ramification locus of $\alpha$. 
Let $L/F(C),\sigma$ and $L'/F(C'),\sigma'$ be the associated ramification 
data. Let $P_c$ be a preimage on $S'$ of $P$, 
which is locally the intersection of $C_c$ and $C_c'$ 
which are preimages of $C$, $C'$ respectively. 
Set $R = {\cal O}_{S,P}$, $R'$ the ring gotten by adjoining 
a primitive $q$ root of one to $R$, and $R'_c$ the localization 
of $R'$ at $P_c$. 
Let $\pi = 0$, $\delta = 0$ define $C$, $C'$ at $R$ 
and similarly for $\pi_c$, $\delta_c$, $C_c$, $C_c'$ and $R'_c$. 
Suppose $M = K((w\pi^s\delta^t)^{1/q})$ where $s$, $t$ 
are prime to $q$. Since $M$ splits the ramification 
of $\alpha$ at $C$ and $C'$, we can define the 
residual Brauer classes $\beta_C \in \Br(F(C))$ 
and $\beta_{C'} \in \Br(F(C'))$ with respect to $M/K$. 

We are interested in describing $M' = M \otimes_K K'$ in terms of $P_c$.  
Since $\pi_c$ appears to the first power in the $R'$ 
factorization of $\pi$, and similarly for $\delta_c$, 
we can write $M' = K'((w_c\pi_c^s\delta_c^t)^{1/q})$. 
Thus there are well defined residual Brauer classes 
$\beta_{C_c} \in \Br(F(C_c))$, $\beta_{C_c'} \in \Br(F(C_c'))$ 
of $\alpha' = \alpha \otimes_K K'$ at $C_c$ and $C_c'$ 
with respect to $M'$. The following is clear. 

\proclaim Proposition 3.10. a) Under the natural maps 
induced by $F(C) \subset F(C_c)$ and $F(C') \subset F(C_c')$, 
$\beta_C$ maps to $\beta_{C_c}$ and $\beta_{C'}$ maps to 
$\beta_{C_c'}$. 
\smallskip
b) Suppose $P$ is a chilly point.  
Then $\beta_C$ and $\beta_{C'}$ are both unramified at $P$ 
and have equal images in $\Br(F(P))$. 
\smallskip
c) Suppose $P$ is a cold point.   
The ramification of $s\beta_C$ and $-t\beta_{C'}$ 
are equal at $P$. $M$ splits all the ramification 
of $\alpha$ at $P$ if and only if the ramification 
of $\beta_C$ is trivial at $P$. 
\smallskip
d) If $\alpha$ has a hot point, 
then the residual classes of $\alpha$ are not split by the ramification. 
In particular, $\alpha$ has index greater than $q$. 

Just as above, we can trivially extend 2.14 as follows.

\proclaim Proposition 3.11. 
Suppose $P$ is a curve point on $C$ 
and the ramification $L/F(C)$ splits at $P$. 
Suppose $M = K(\pi^{1/q})$ and $\pi$ has $C$ valuation 
prime to $q$. If $\beta_C$ is the residual Brauer class at $C$ with respect 
to $M$, then $\beta_C$ is unramified at $P$.

\proof 
This is obvious by functoriality, 2.8, 2.12 and 2.13. 

\bigskip
\leftline{Section Four: Killing the residual class} 
\medskip
In Proposition 0.6 we saw how one can modify the residual 
class by changing the ramified extension. Next we observe 
how we can do that for several curves at once. 
To this end, let $S$ be an excellent nonsingular surface projective over some 
affine $A$. Set $K = F(S)$ and suppose $\alpha \in \Br(K)$ 
is of order $q$. 

We need to be slightly more general about the ramification locus. 
Let $B$ be a finite set of curves on $S$ including the ramification 
locus. 
As usual, suppose we have blown up $S$ so that 
$B$ consists of smooth curves with 
normal crossings. Let $\{P_j\}$ be the set of nodal 
points on the ramification locus, and assume 
we have further blown up so that there are no chilly 
loops and no cool points. Set $\{L_i/F(C_i),\sigma_i\}$ 
to be the ramification data of $\alpha$. 
Suppose that for each 
$C_i$ in the genuine ramification locus 
we fix $s_i$, as in 3.8, such that $s_i$ is prime to $q$ 
and the following holds. If $P_j$ is a chilly point which is locally 
the intersection of $C_i$ and $C_j$, and $s$ is the 
coefficient of $\alpha$ at $P$ with respect to $C_i$, 
then $s = s_{i'}/s_i$ in $\Z/q\Z$. 

Let ${\cal P}$ be a finite set of closed points 
including all nodal points of $B$. If any curve of $B$
contains only finitely many closed points, we can assume 
${\cal P}$ contains them all. 

\proclaim Lemma 4.1. Let ${\cal P}$ be as above. 
We can choose $\pi \in K$ such that the support of 
$E = (\pi) - \sum_i s_iC_i$ contains no components of $B$, 
only intersects $B$ in nonsingular points, and contains no 
point of ${\cal P}$. 

\proof 
Use weak approximation to choose $\pi'$ with valuation 
$s_i$ at $C_i$. Write $(\pi) = \sum_i s_i C_i + E$. 
We can assume ${\cal P}$ includes a point on every component of $B$. 
By 1.5 there is a $u \in K$ with $(u) = E' - E$ 
where the support of $E'$ does not contain any element 
of ${\cal P}$. Now $\pi = u\pi'$ is as needed.~\qed 

Let $s_i$ and $\pi$ be as in 4.1. 
Set $M = K(\pi^{1/q})$.  
Let $\beta_{C_i}$ be the residual Brauer classes at $C_i$ 
with respect to $M$. In the rest of this section we assume 
all the residual Brauer classes of $\alpha$ at the $C_i$ 
are split by the ramification. By 0.5 this happens 
if $\alpha$ has index $q$. 
Note that this assumption means 
$\beta_{C_i} = \Delta(L_i/F(C_i),\sigma_i,u_i)$ 
for some $u_i \in F(C_i)^*$. 

First we consider the ramification of the $\beta_{C_i}$ 
at non-nodal points.

\proclaim Theorem 4.2. Let $\pi$ be as above, and 
$C_i$ some curve in the ramification locus of $\alpha$. 
Let $\beta_{C_i} = \Delta(L_i/F(C_i)$, $\sigma_{C_i}, u_i)$  
be as above. Let $P$ be a non-nodal point on $C_i$. 
\smallskip
a) If $P$ is not in the support of $E$,  
then $\beta_{C_i}$ has is unramified at $P$. 
\smallskip
b) Suppose $P$ is in the support of $E$. Then 
$L_i/F(C_i)$ is unramified at $P$. If $L_i/F(C_i)$ is split 
at $P$, then $\beta_{C_i}$ is again unramified at $P$. 
\smallskip
c) Suppose $P$ is in the support of $E$ and $L_i/F(C_i)$ 
is not split at $P$. Let $\gamma = \bar L_i/F(P), \bar \sigma_{C_i}$ 
be the induced extension of $F(P)$ viewed as an element of 
$H^1(F(P), \Q/\Z)$. Then the ramification of $\beta_{C_i}$ 
has the form $-m_i(C_i.E)_P\gamma$ 
where $(C_i.E)_P$ is the intersection multiplicity at $P$ 
and $m_i$ is the modulo $q$ inverse of $s_i$. 

\proof 
Let $R = {\cal O}_{S,P}$. By the usual trick, it suffices to prove this theorem 
after adjoining a primitive $q$ root of one, $\rho$, which we fix. 
Let $\pi_i \in R$ be a prime of $R$ defining $C_i$ 
locally at $P$. We know that $\alpha = \alpha' + (u,\pi_i)$ 
for some $u \in R^*$ and $\alpha' \in \Br(R)$. 
Then if $\bar u$ is the image of $u$ in $F(P)$, 
$\bar L_i/F(P), \bar {\sigma_i}$ is the same as 
$F(P)(\bar u^{1/q})/F(P),\bar {\sigma_i}$ 
where $\bar {\sigma_i}(\bar u^{1/q})/\bar u^{1/q} = \rho$.

We turn to proving a). 
Perhaps up to $q$ powers, 
$\pi = v\pi_i^{s_i}$ where $v \in R^*$. 
It follows that for some $u' \in R^*$, 
$\alpha = \alpha'' + (u',\pi)$ where $\alpha'' \in \Br(R)$.  The elements 
$\alpha$ and and $\alpha''$ have the same image in $\Br(M)$ and  
we can use $\alpha''$ to compute $\beta_{C_i}$. Since $\alpha'' \in \Br(R)$, 
$\beta_{C_i}$ is unramified at $P$. 

Next we prove b). 
Set $E_P$ to be the sum  $\sum t_jE_j$ over all 
$E_j$ in the support of $E$  which intersect 
$C_i$ at $P$. For each $E_j$ in the support of $E_P$ 
let $\delta_j \in R$ be a prime 
such that $\delta_j = 0$ 
defines $E_j$ at $P$. Set $\delta = \prod \delta_j^{t_j}$, 
the product over the support of $E_P$.  
Then, up to $q$ powers, $\pi = v\pi_i^{s_i}{\delta}$ where $v \in R^*$.  

Let $s_im_i - 1$ be divisible by $q$, so up to $q$ powers 
$\pi_i$ is $\pi^{m_i}(v\delta)^{-m_i}$. 
The element $\alpha$ can be rewritten as $\alpha' + 
(u^,(v\delta)^{-m_i}) + 
(u^{m_i},\pi)$. 
As before, $\alpha$ has the same image in $\Br(M)$ as 
$\alpha' + (u,v^{-t_i'}\delta^{-t_i'})$ and the 
image of $\alpha'$ is unramified at $P$. If $L_i/F(C_i)$ 
is split at $P$,  
then $\bar u$ is a $q$ power in $F(P)^*$ and 
$\beta_{C_i}$ again is unramified at $P$. 
This proves b). 
Otherwise by 0.12, the ramification of 
$\beta_{C_i}$ is defined by $(\bar u^{-m_in})^{1/q}$ where $n$ 
is the valuation of $\bar \delta$ at $P$, and hence is 
$(C_i.E)_P$.~\qed  

We fix ${\cal Q}$ to be a finite set of closed points 
on the ramification locus 
which are on only one $C_i$ and where the relevant 
$\beta_{C_i}$ are unramified. If $P$ is a point on 
a $C_i$ and a component of $B$ not among the $C_i$, 
then by 4.1 and 4.2 the relevant $\beta_{C_i}$ 
is unramified at $P$ and we can assume $P$ is in ${\cal Q}$. Furthermore, by 4.1 and 4.2 we can assume 
that any curve among the $C_i$ with no nodal points at all 
contains a point of ${\cal Q}$.

\proclaim Proposition 4.3. Let ${\cal Q}$ be a finite 
set of closed points as above.  
Assume all the residual Brauer classes of $\alpha$, the  
$\beta_{C_i}$, are split by the ramification, so $\beta_{C_i} = 
\Delta(L_i/F(C_i),\sigma_i,u_i)$.  In particular, assume
there are no hot points. Let 
Then there are $v_i \in F(C_i)$ such that:
\smallskip
i) The $v_i$ are units at all nodal points and all the $Q_l$. 
\smallskip
ii) $\Delta(L_i/F(C_i),\sigma_i,v_i) = \Delta(L_i/F(C_i),\sigma_i, u_i^{s_i})$.
\smallskip
iii) If $P$ is a nodal point and at the intersection of $C_i$ 
and $C_{i'}$, then $v_i$ and $v_{i'}$ have equal images 
in $F(P)$. 

\proof Let $P_j$ be a nodal point on $C_i$. 
Let $\hat F_j$ be the completion of $F(C_i)$ at $P_j$ 
and let $\hat L_j = L_i \otimes_{F(C_i)} \hat F_j$. 
Define $N_i$ to be the norm map of $L_i/F(C_i)$ and 
$\hat L_j/\hat F_j$. 
If $\hat L_j$ is split or ramified at $P_j$, 
then norms have all possible valuations, 
so we can choose $w_j \in \hat L_j$ 
such that $N_i(w_j)/u_i$ is a unit. If $\hat L_j$ is a field 
and unramified at $P$, then $P$ must be a chilly point 
and $\Delta(L_i/F(C_i),\sigma_i,u_i)$ is unramified 
at $P$ (3.10). Thus  $u_i$ 
must have valuation a multiple of $q$ and we can choose 
$w_j$ such that $N_i(w_j)/u_i$ is a unit. 
At a point of ${\cal Q}$ we have assumed $\Delta(L_i/F(C_i),\sigma_i,u_i)$ 
is unramified and so once again $w_l$ exists with $N_i(w_l)/u_i$ a unit at that point. 
By weak approximation 
we can find $w \in L_i$ such that $N_i(w)/N_i(w_j)$ is a unit at all 
nodal points $P_j$ and all points of ${\cal Q}$. 
Then $u_i$ can be replaced by 
$u_i/N_i(w)$ and we can assume all the $u_i$ are units at 
all the nodal points and all the points of ${\cal Q}$. 

For clarity's sake, set $v_i' = u_i^{s_i}$. 
Let $v_i'(P)$ be the image of $v_i'$ in the residue field 
$F(P)$ of a point $P$ on $C_i$ (when defined). Suppose $P_j$ is a nodal chilly 
point at the intersection of $C_i$ and $C_{i'}$ with 
coefficient $s$ with respect to $C_i$.  
Let $\bar L_{ij}/F(P_j), \sigma_{ij}$ be the residue extension of 
$L_i/F(C_i), \sigma_i$ 
at $P_j$. This is well defined, a field, and equal to 
$\bar L_{i'j}/F(P_j),\sigma_{ij'}^s$, by the definition of $s$ and 
chilly point. 

\proclaim Lemma 4.4. $\,\,$ 
If $P_j$ is a chilly point, 
$v'_i(P_j)$ and $v'_{i'}(P_j)$ differ by a norm 
of $\bar L_{ij}/F(P) = \bar L_{i'j}/F(P)$. If 
$P_j$ is a cold point, $v'_i(P_j)$ and $v'_{i'}(P_j)$ 
differ by a $q$ power. 

\proof 
If $P_j$ is a chilly point, we know by 3.10 and 3.9  that 
$$\Delta(\bar L_{ij}/F(P_j),\sigma_{ij},v'_i(P_j)) = 
s_i\Delta(\bar L_{ij}/F(P_j),\sigma_{ij},u_i(P_j))$$ 
equals 
$$s_i\Delta(\bar L_{i'j}/F(P_j),\sigma_{i'j},u_{i'}(P_j)) = 
s_{i'}s\Delta(\bar L_{ij}/F(P_j),\sigma_{i'j},u_{i'}(P_j)) = $$
$$ = \Delta(\bar L_{ij}/F(P_j),\sigma_{ij},v'_{i'}(P_j))$$ 
because 
$\sigma_{ij'}^s = \sigma_{ij}$ and $v'_{i'} = u_{i'}^{s_{i'}}$. 
By e.g. [LN] p. 45, we are done for chilly $P_j$. 

If $P_j$ is a cold point, we know by 3.10 that 
$$s_i\Delta(L_i/F(C_i),\sigma_i,u_i)\ \hbox{ and }\ 
s_{i'}\Delta(L_{i'}/F(C_{i'}),\sigma_{i'},u_{i'})$$ 
have inverse ramifications at $P_j$. Moreover, by 
3.7, $L_i/F(C_i),\sigma_i$ and $L_{i'}/F(C_{i'}),\sigma_{i'}$ 
have inverse ramifications at $P_j$. It follows from 0.10 
that $v'_i(P_j)$ and $v'_{i'}(P_j)$ differ by a $q$ 
power.~\qed 

We are ready to finish the proof of 4.3, which is now easy. 
Weak approximation implies that we can modify the $v'_i(P)$ 
by norms or $q$ powers independently at all the nodal points and a
all points in ${\cal Q}$. Proposition 4.3 is immediate.~\qed 

The point of 4.3 was to have enough compatibility 
among the $v_i$ to do the following. 

\proclaim Proposition 4.5. Let $v_i$ and $Q_l$ be as in 4.3, 
and continue the same assumptions on the residual Brauer classes 
and the lack of hot points.   
Then there is an affine $U \subset S$ with 
affine ring $R$ and a $v \in R$ such that the 
following holds. 
\smallskip
a) $U$ contains all nodal points, contains ${\cal Q}$, and contains a closed point on all the curves in $B$. 
\smallskip
b) If $R_i$ is the affine ring of $U \cap C_i$, then 
$v_i \in R_i$. 
\smallskip
c) The element $v$ is a unit at all curves of $B$, at 
all nodal points of $B$ and maps to $v_i$ for all 
$i$. 

\proof We can choose a set ${\cal P}$ of closed points so that 
the following is true. First, the points of ${\cal P}$ are not on any 
$C_i$, have a point on any component of $B$ not intersecting a 
$C_i$, and include all nodal points of $B$ not on any $C_i$. 
Thus among the points of ${\cal P}$, ${\cal Q}$, and the nodal points of the 
ramification locus are all nodal points of $B$ and at least 
one point on any component of $B$. 

By 1.4 there is an affine open $U' \subset S$ 
containing ${\cal P} \cup {\cal Q}$, 
and containing all the nodal points of the $C_i$. 
Let $P'_n$ be the set of poles of the $v_i$ on $U \cap C_i$. 
There is an $f$ defined on $U$ which is 0 on all the 
$P'_n$ and nonzero at all $P_m$, all 
the nodal points, and all the $Q_l$. We set $U = U'_f$. 
This finishes a) and b). 

Let $Q_i \subset R$ be prime ideals corresponding 
to the $C_i$ and the $P_m$. If $Q_i$ corresponds to a 
$P_m$, let $v_i$ be arbitrary nonzero. 
Note that the $Q_i$ corresponding to a $P_m$ are maximal and 
relatively prime to any other of the $Q_i$. 
Translating in commutative algebra, 
we have a ring $R$, prime ideals $Q_i$ with no 
inclusions among them, and 
elements $v_i \in R/Q_i$ such that the following holds. 
First, $Q_{i} + Q_{i'}$ is either $R$, or a finite 
intersection of maximal ideals $M_j$ and each 
maximal ideal contains at most two $Q_i$. Second, whenever 
$M_j$ contains $Q_i + Q_{i'}$, $v_i$ and $v_{i'}$ have equal 
images in $R/M_j$. Just using these facts, c) is proven by induction 
on the cardinality of the set of $Q_i$. Of course one $Q_i$ is trivial. 
Suppose $v'$ is chosen for $Q_1,\ldots,Q_{n-1}$. 
Set $J = \cap_{i=1}^{n-1} Q_i$ and $I = Q_n$. 
We claim $I + J$ is the intersection of maximal ideals 
$M_j$ where $M_j$ contains $Q_n$ and one of the $Q_i$, $i < n$. 
But $I + J \subset \cap_j M_j$ is clear, and equality can be shown 
by checking it locally. But $R/(I + J)$ is the direct sum 
of the $R/M_j$ and c) follows from the exact sequence 
$0 \to R/(I \cap J) \to R/I \oplus R/J \to R/(I + J) \to 0$.~\qed 

\proclaim Theorem 4.6. Let $\alpha$, $C_i$, $Q_l$, and 
$s_i$ be as above.  Assume all the residual Brauer classes 
at all the $C_i$ are split by the ramification, and hence that
there are no hot points. 
Then there is a new choice of $\pi \in K$, 
such that $\pi$ has valuation $s_i$ at the $C_i$, 
$E = (\pi) - \sum_i s_iC_i$ does not contain any nodal 
points of $B$,  or any point in ${\cal Q}$, 
or any components of $B$  
in its support, and with respect to $M = K(\pi^{1/q})$, 
all of the residual Brauer classes $\beta_{C_i}$ 
are trivial. 
Furthermore, $(C_i.E)_P$ is a multiple of $q$ for all points $P$ 
on the $C_i$ where $L_i/k(C_i)$ is nonsplit. 

\proof We find $\pi'$ as in 4.1 and $v$ as in 4.5. 
Then by 0.7 $\pi = v\pi'$ has all the residue classes split. 
The last sentence follows from 4.2 c).~\qed 

Combining 4.6 with 3.9 and 3.10 we have: 

\proclaim Corollary 4.7. If $M$ is as in 4.6, 
then $M$ splits all the ramification of $\alpha$ 
at chilly and cold points. 

\bigskip
\leftline{Section 5: The proof}
\medskip
We have gone as far as we can assuming $S$ is a fairly general 
surface. In this section, we return to the situation 
of section one and assume $S \to \Spec(\Z_p)$ is projective, regular, 
excellent, of finite type, with relative dimension one. 
Let $\bar S$ be the reduced subscheme defined by the 
preimage of the closed point of $\Spec(\Z_p)$.  
We assume $\alpha \in \Br(K(S))$ has order $q$, 
and we let $B$ be the union of the ramification 
locus of $\alpha$ and $\bar S$. 
We can blow up $S$ so that 
$B$ consists of nonsingular curves with normal crossings, and 
so that there are no cool points or chilly loops. 
If $C$ is any curve on $S$, 
then the residue field of $C$ will be written as 
$k(C)$ to emphasize that it is a curve over a finite 
field or $\Spec$ of a $p$ adic number ring. 
For each $C_i$ in the ramification locus of $\alpha$, let 
$L_i/k(C_i), \sigma_{C_i}$ be the ramification data. 
Until 5.2 we assume that all the residual Brauer classes of $\alpha$ 
are split by the ramification, and hence that there are no hot points. 

Let $\pi$ be as in 4.6 and write (again) $(\pi) = \sum s_iC_i + E$.  
Let $\bar E$ be the divisor which is the sum, with 
coefficients 1, of all the curves in $\bar S$. 
Let $\gamma \in \Pic(S)$ be the line bundle equivalent 
to the divisor class $-E$, and 
$\bar \gamma \in \Pic(\bar S)$ its image. 
Then $E$ and $\bar E$ only intersect 
in smooth points of $\bar S$ and so we can represent $\bar \gamma$ 
as a divisor using the intersection of $-E$ and $\bar E$. 
In particular, $\bar \gamma$ has the form $\sum_j qn_jQ_j + \sum_l n_lQ_l'$ 
where by 4.6 the $Q_l'$ are either not on the ramification locus of $\alpha$ 
or are at points where $L_i/k(C_i)$ splits. For each of the 
$Q_l'$ choose a geometric curve $E_l' \subset S$ whose unique 
closed point is $Q_l'$ (1.1). Set $E' = -E - \sum_l n_lQ_l'$. In the notation 
of 1.6, let $P$ represent the set of all nodal points 
on $B$. 
Consider the element $\gamma' \in  H^1(S,{\cal O}^*_P)$ represented, 
as in 1.6, by the divisor $E'$ and the element $1$ at all points in $P$.  
 
The image, $\bar \gamma'$, of $\gamma'$ in $H^1(\bar S,{\cal O}^*_P)$ 
lies in $qH^1(\bar S,{\cal O}^*_P)$. It follows 
from 1.7 that $\gamma$ lies in $qH^1(S,{\cal O}^*_P)$. 
That is, using 1.6, there is a divisor $E''$, elements $a_j \in k(P_j)^*$ 
for all $P_j$ in $P$, and an $f \in K = F(S)$ such that $f$ is a unit at all 
nodal points, $(f) = E' + qE''$ and $f(P_j) = a_j^q$ at all $P_j$. 

Now we compute the divisor $(f\pi) = \sum_i s_iC_i + \sum n_jD_j$. 
We note that for any curve $D_j$, $D_j$ intersects $B$ in a smooth point, 
and if $n_j$ is prime to $q$, $D_j$ either does not intersect 
any $C_i$, or does so at a point where $L_i/F(C_i)$ splits. 

\proclaim Theorem 5.1. Let $K$ be a field finite over $\Q_p(t)$. 
Let $\alpha \in \Br(K)$ have index a prime $q \not= p$. 
Then $\alpha$ is represented by a cyclic algebra of degree $q$. 

\proof As in [S], we know $K$ is the function field 
of a regular excellent projective surface $S$ projective 
over $\Spec(\Z_p)$. As we have said before, we can blow up 
so that $B$, the union of the ramification locus of $\alpha$ 
and $\bar S$, consists of regular curves with normal 
crossings. We can further blow up so that the ramification locus has no cool 
points or chilly loops. By the assumption on the 
index, there are no hot points and the residual classes 
are all split by the ramification. 
Find $\pi$ as 4.6. Choose $f$ as above, 
and write $M = K((f\pi)^{1/q})$. 
For each curve $C_i$ in the ramification locus, 
let $\beta_{C_i}$ be the residual Brauer class of 
$\alpha$ at $C_i$ with respect to $M/K$. 
We claim 
$\alpha' = \alpha \otimes_K M$ is not ramified on any 
discrete valuation over $S$. 

The choice of $s_i$ insures that $\alpha'$ is not ramified 
on the primes over the $C_i$, the curves in 
the ramification locus of $\alpha$. Since 
$\alpha$ itself is unramified at all other curves, 
we are reduced to considering discrete valuations 
over points of $S$. By 3.4 
we can also ignore distant points and curve points 
$P \in C_i$ where the ramification $L_i/F(C_i)$ 
splits.  
If $M' = K(\pi^{1/q})$, then by 4.6 all the residual 
classes with respect to $M'/K$ are trivial. 
Since $f(P_j) \in (k(P_j)^*)^q$, it follows from 
0.7, 3.9, and 3.10 that $\alpha'$ is unramified at any discrete 
valuation over a nodal point. 
Finally suppose 
$P$ is a curve point on $C_i$ where the ramification is nonsplit. 
By our choice of $f\pi$, the only curves in the support 
of $(f\pi)$ that meet $P$ have coefficients a multiple of $q$. 
That is, if $R = {\cal O}_{S,P}$, then $f\pi = u\pi_c^s\delta^q$ 
where $u \in R^*$, $\pi_C = 0$ defines $C$ locally at $P$, 
and $s$ is prime to $q$. By 3.5, $M$ splits all the ramification. 
By 0.9, $M$ splits $\alpha$, and so by 0.1 $\alpha$ is represented 
by a cyclic algebra of degree $q$.~\qed

One might be interested in how to detect those $\alpha$ of index $q$. 
The answer is not complicated. 

\proclaim Corollary 5.2. Suppose $S$ is as in this section, 
$K = F(S)$, 
and $\alpha \in \Br(K)$ has order $q$ in the Brauer 
group. Assume $S$ has been blown up so that 
the ramification locus of $\alpha$ 
consists of nonsingular curves with normal crossings. 
Then $\alpha$ has index $q$ if and only if 
there are no hot points. 

\proof 
Up until the statement of 5.1 we only assumed that all the residual 
Brauer classes were split by the ramification. We did not make this 
part of 5.1 only because it would be clumsy to state. 
So to prove 5.2, it suffices to show that without hot points, 
all the residual Brauer classes are split by the ramification. 
Consider $C$ in the ramification locus, and let 
$M = K(\pi^{1/q})$ where the $C$ defined valuation of $\pi$ 
is prime to $q$. Set $\beta_C$ to be the residual Brauer class 
with respect to $M$, and let $L/k(C),\sigma$ be the ramification 
of $\alpha$ at $C$.  Since $\beta_C$ must have order $q$ 
(or 1), by 0.8 to show $L$ splits $\beta_C$ it suffices to show 
$L$ splits the residues of $\beta_C$ at all points $P$. This is 
automatic at any point where the prime defining $P$ 
does not split in $L$ (e.g. use 0.3). Thus it suffices 
to show $\beta_C$ is unramified at all points where $L/k(C)$ 
splits. But this is 3.11.~\qed

\def\hangbox to #1 #2{\vskip1pt\hangindent #1\noindent \hbox to #1{#2}$\!\!$} 
\def\ref#1{\hangbox to 30pt {#1\hfill}}

\leftline{References}
\bigskip
\ref{[A]} Albert, A. Adrian, ``Structure of Algebras'', 
American Mathematical Society, Providence RI 1961 (Colloquium Publ. v. 24)
\bigskip
\ref{[E]} Eisenbud, David, ``Commutative Algebra with a view toward Algebraic Geometry'', Springer-Verlag New York 1995 
\bigskip
\ref{[EGA]} Grothendieck, A. and Dieudonne, J., 
``El\'ements de G\'eom\'etrie Alg\'ebrique III'' 
{\it Etude cohomologique des faisceaux coh\'erent}, 
Publ. Math. IHES 11 (1961)
\bigskip
\ref{[G]} Grothendieck, A.,
{\it Le groupe de Brauer III: exemples et compl\'ements},
in ``Dix Expos\' sur la Cohomologie des sch\'emas'', 
North Holland, Amsterdam, 1968 
\bigskip
\ref{[H]} Hartshorne, Robin, ``Algebraic Geometry'', 
Springer-Verlag, New York/\hfill\break Heidelberg/Berlin, 1977 
\bigskip
\ref{[H1]} Hartshorne, Robin, ``Ample Subvarieties of Algebraic Varieties'', 
Springer-Verlag, Heidelberg, 1970 (Lecture Notes in Mathematics 156)
\bigskip
\ref{[JW]} Jacob, Bill and Wadsworth, Adrian, 
{\it Division Algebras over Henselian Fields}, 
J. of Algebra {\bf 128} no. 1  (1990) 126--179 
\bigskip 
\ref{[L]} Lipman, J.,
{\it Introduction to resolution of singularities},
in Hartshorne, R. ed., ``Algebraic Geometry Arcata 1974'' 
American Mathematical Society, Providence RI, 1975
(Proceedings in Pure Mathematics v. 29 p. 187--230)
\bigskip
\ref{[LN]} Saltman, David J., ``Lectures on Division Algebras'', 
American Mathematical Society, Providence, RI, 1999 
(CBMS Lecture Note Series \#94) 
\bigskip
\ref{[Mi]} Milne, J.S., ``Etale Cohomology'', Princeton University 
Press, Princeton, N.J., 1980
\bigskip
\ref{[R]} Reiner, Irving, ``Maximal Orders'', Academic Press, 
London/New York/San Francisco, 1975
\bigskip
\ref{[S]}
Saltman, David J., {\it Division algebras over p-adic curves}, 
J. Ramanujan Math. Soc., {\bf 12}  (1997), 25--47 
and {\it Correction to ``Division algebras over p-adic curves''}, 
J. Ramanujan Math. Soc. {\bf 13} (1998), 125--129 
\bigskip
\ref{[Se]}
Serre, J.P.,  ``Local Fields'', 
Springer-Verlag, New York/Heidelberg/Berlin, 1979  

\end